\documentclass[12pt]{article}
\usepackage{amsmath,amssymb}
\newcommand{\qed}{{\unskip\nobreak\hfil\penalty50\hskip2em\vadjust{}
    \nobreak\hfil$\square$\parfillskip=0pt\finalhyphendemerits=0\par}}

\newcommand{\x}{\times}
\newcommand{\<}{\langle}
\renewcommand{\>}{\rangle}
\renewcommand{\a}{\alpha}
\renewcommand{\b}{\beta}
\newcommand{\bs}{\bigskip}
\newcommand{\dt}{\cdot}
\newcommand{\C}{\mathbb C}
\newcommand{\Cn}{\mathbb C_{\langle n\rangle}}
\newcommand{\ds}{\displaystyle}
\renewcommand{\d}{\delta}
\newcommand{\D}{\Delta}
\newcommand{\e}{\varepsilon}
\newcommand{\g}{\gamma}
\newcommand{\G}{\mathcal G}
\newcommand{\F}{\mathcal F}
\renewcommand{\i}{\infty}
\renewcommand{\l}{\lambda}
\renewcommand{\L}{\Lambda}
\newcommand{\LL}{\mathcal L}
\newcommand{\mb}{\mbox}

\newcommand{\N}{\mathbb N}
\renewcommand{\O}{\Omega}
\renewcommand{\o}{\omega}
\newcommand{\ot}{\otimes}
\newcommand{\p}{\partial}

\newcommand{\stt}{\subset}
\newcommand{\stk}{\stackrel}
\newcommand{\T}{\mathcal T}
\renewcommand{\th}{\theta}
\newcommand{\var}{\varphi}
\newcommand{{\z}}{\Bbb Z}
\newcommand{\VE}{\mbox{Vect}\,{\mathcal E}}
\newcommand{\R}{\mathbb R}
\newcommand{\alg}{algebra}
\newcommand{\pr}{probabili}
\newcommand{\po}{polynomial}
\newcommand{\ncc}{noncommutative }

\newcommand{\dis}{distribution}
\newcommand{\gra}{gradient}
\newcommand{\mor }{morphism}
\begin{document}
\setlength{\baselineskip}{16pt}

\begin{center}
{\Huge Cyclomorphy}

\bs\bs {\sc Dan Voiculescu\footnote{Research at
MSRI supported in part by
National Science Foundation grant DMS--0079945.\newline}}

\bs Department of Mathematics\\University of California\\
Berkeley, California 94720--3840
\end{center}

\bs
\begin{quote}
{\bf Abstract.} A main result is that, roughly, a dense set of the
infinitesimal trace-preserving deformations of a semicircular system
$s_1,\dots ,s_n$ arise from one-parameter groups of auto\mor s of
the free-group factor $L(F(n))$ generated by $s_1,\dots ,s_n$.
More generally the paper studies cyclic \gra s in von Neumann \alg s,
Lie \alg s of \ncc trace-preserving vector fields and the class of
cyclomorphic maps which preserve the orthogonals of spaces of cyclic
\gra s.\end{quote}
\bs

\setcounter{section}{-1}
\section{Introduction}
The word ``cyclomorphy" is meant to remind the reader of
``holomorphy" and of the ``cyclic derivative" of Rota-Sagan-Stein ([8]).
The paper is about cyclic derivatives in von Neumann \alg s.

Several questions around free \pr ty theory lead to cyclomorphy:
moments of \ncc random variables, free entropy ([12],[14]) and its connection
with large deviations for random matrices ([4],[5]), orbits of the equivalence
in \dis \ and the free Wasserstein distance ([1]). Last but not least I was
motivated by curiositiy about Aut$(L(F(n)))$, the auto\mor \ group of a free
group factor.

In general, if $Y_j=Y^*_j$, $1\leq j\leq m$, generate a II${}_1$ factor $M$,
the orbit of $Y=(Y_j)_{1\leq j\leq n}$ under Aut$(M)$ provides a
parametrization of Aut$(M)$. The derivations in the ``Lie algebra of Aut$(M)$"
which have the $Y_j$'s in their domain of definition, are then parametrized
by the ``tangent space" of the orbit at $Y$. This ``tangent space" is a more
tractable object as we shall see.

Since auto\mor s of $M$ preserve its trace-state $\tau$, the``tangent space" is
orthogonal to the space of cyclic gradients of \ncc \po s in
$(Y_1,\dots ,Y_n)$. The question then is which of the self-adjoint $n$-tuples
in the orthogonal of cyclic \gra s give rise to derivations which exponentiate
to one-parameter groups of auto\mor s? We show that exponentiation is possible
if the $n$-tuple is \po , or more generally is given by \ncc power series of
sufficiently large radius of convergence evaluated at $Y$.

The \po \ elements in the orthogonal to cyclic \gra s form a Lie \alg , which
is a natural \ncc relative of the Lie \alg \ of vector fields which preserve a
volume form.

  A natural generator for $L(F(n))$ is provided by a semicircular system
$s_1,\dots ,s_n$ ([10],[11]). In  this case we can give a complete description
of the orthogonal to cyclic \gra s. In particular we show that \po \ elements
are dense in this space. Exponentiation therefore works for the corresponding
derivations and we have an ``infinitesimally rich" sub\alg \ of the ``Lie \alg
" of the auto\mor \ group in this case.

 From this study of automorphic orbits it becomes clear that it is important to
understand the differential geometric picture arising on $M^n$ (or on the
hermitian subspace $M^n_h$) from the field of subspaces of the cotangent spaces
arising from cyclic \gra s and their orthogonals which are a field of subspaces
in the tangent spaces. In particular there is a natural class of maps, those
for which the differentials preserve the orthogonals to cyclic \gra s.
We look at basic properties of cyclomorphic maps and of a generalization of
these to the case of $B$-von Neumann \alg s ($B$ a von Neumann sub\alg ).

The paper has eight sections without the introduction. Section 1 deals with
preliminaries. Section 2 is about the estimates which will prove that certain
derivations can be exponentiated to auto\mor \ groups.
Section 3 combines the estimates of section 2 with the condition of
orthogonality to cyclic \gra s to obtain results about exponentiation of
derivations.
Section 4 is a brief look at endomorphic orbits and their tangent spaces.
In section 5 we study real and complex cyclomorphic maps. In Theorem 5.13
and Proposition 5.14 we point out the role of independence and of free
independence in this context.
Section 6 is a collection of generalities about the Lie \alg s of ``\ncc
vector-fields" which occur in this paper. Section 7 deals with the semicircular
case. We compute the orthogonal to the cyclic \gra s in this case.
Section 8 is about $B$-morphic maps, the extension to von Neumann \alg s
``over $B$" of cyclomorphic maps. We only sketch how this generalization is
done and will develop the details elsewhere.

\section{Preliminaries}
\subsection{Noncommutative polynomials}

If $X_1,\dots ,X_n$ are noncommuting indeterminates, the ring of \ncc \po s
will be denoted by ${\C}_{\<n\>}={\C}\<X_1,\dots ,X_n\>$.
If $a_1,\dots, a_n$ are elements in an unital \alg \ over ${\C}$ and
$P\in{\C}_{\<n\>}$ then $P(a_1,\dots ,a_n)$ or
$\e_{(a_1,\dots ,a_n)}P$ denotes the evaluation of $P$ at $a_1,\dots ,a_n$.

We shall also use the notation ${\C}\<a_1,\dots ,a_n\>$ for
$\e_{(a_1,\dots ,a_n)}{\C}\<X_1,\dots ,X_n\>$, i.e.~for the \alg \ generated
by $1,a_1,\dots ,a_n$. Clearly this \alg \ is isomorphic to the
\alg \ of \ncc \po s only if $a_1,\dots ,a_n$ are \alg ically free.
When we will use the notation it will be clear from the context whether
$a_1,\dots ,a_n$ are \alg ically free or not.

\subsection{Cyclic gradients and free difference quotients}

The partial free difference quotients are the derivations
\[
\p_j:{\C}_{\<n\>}\to {\C}_{\<n\>}\ot {\C}_{\<n\>}
\]
such that $\p_jX_k=0$ if $j\neq k$ and $\p_jX_j=1\ot 1$,
where $j,k\in\{1,\dots ,n\}$. The $\Cn$-bimodule structure on
${\Cn}\ot{\Cn}$ is the obvious one $a_1(b\ot c)a_2=a_1b\ot ca_2$.

The partial cyclic derivatives are then
$\d_j ={\tilde\mu}\circ\p_j:{\Cn}\to{\Cn}$ where
${\tilde\mu}(a\ot b)=ba$. Thus
\[
\d_j X_{i_1}\dots X_{i_p} =\sum_{1k|1\leq k\leq p, \ i_k=j}
X_{i_{k+1}}\dots X_{i_p}X_{i_1}\dots X_{i_{k-1}} \ .
\]
The map $\d:{\Cn}\to ({\Cn})^n={\Cn}\oplus\dots\oplus{\Cn}$ given by
\[
\d P=(\d_j P)_{1\leq j\leq n} \ = \ \d_1 P\oplus\dots\oplus \d_n P
\]
is the cyclic gradient.

Similarly $\p:{\Cn}\to ({\Cn}\ot{\Cn})^n$ denotes the free difference
quotient \gra.

In case $a_1,\dots ,a_n$ are \alg ically free the corresponding maps
for ${\Cn}\<a_1,\dots ,a_n\>$ will be denoted
\[
\p_j^{(a_1,\dots ,a_n)},\p^{(a_1,\dots ,a_n)},
\d_j^{(a_1,\dots ,a_n)},\d^{(a_1,\dots ,a_n)}
\]
or $\p^\a_j,\p^\a,\d^\a_j,\d^\a$ where $\a=(a_1,\dots ,a_n)$.

The notation for cyclic derivatives and gradients should not be confused with
Kronecker symbols $\d_{ij}$ which have two lower indices.

If $a,a_1,\dots ,a_n$ are elements in some unital \alg \ $A$ over ${\C}$
let $m_a:{\C}\<a_1,\dots ,a_n\>\ot{\C}\<a_1,\dots ,a_n\>\to A$ be the map
given by
\[
m_a(P_1\ot P_2) \ = \ P_1aP_2 \ .
\]
In particular if $a_1,\dots ,a_n$ are \alg ically free, then a derivation of
${\C}\<a_1,\dots ,a_n\>$ into $A$ can always be written
\[
P\rightsquigarrow \sum_{1\leq j\leq n} m_{b_j} \, \p^\a_j P
\]
where $b_1,\dots ,b_n\in A$. The elements $b_1,\dots ,b_n$ are uniquely
determined, being the values of the derivation on $a_1,\dots ,a_n$.

\subsection{The exact sequence for cyclic \gra s}

In [13] we described the set of cyclic \gra s by an exact sequence
\[
0\to {\C}1+[{\Cn},{\Cn}]\to {\Cn}\stk{\d}{\to} ({\Cn})^n\stk{\th}{\to}
{\Cn}
\]
where $\th(P_1\oplus\dots\oplus P_n)=\sum_j [X_j,P_j]$.
Moreover
\[
{\C}1+[{\Cn},{\Cn}] = {\C}1+\sum_{1\leq k\leq n}
[X_k,{\Cn}] = {\mb{Ker}} \ C
\]
where \ $C:{\Cn}\to{\Cn}$ is the cyclic symmetrization map
$C1=0$ and
\[
C\,X_{i_1}\dots X_{i_p} =\sum_{1\leq j\leq p}
X_{i_{j+1}}\dots X_{i_p}X_{i_1}\dots X_{i_j}
\]
if \ $p\geq 1$.

If instead of $X_1,\dots ,X_n$ we have \alg ically free $a_1,\dots ,a_n$
we shall use the notations
\[
\th^\a \ , \ \ C^\a
\]
where $\a=(a_1,\dots ,a_n)$.

Also, if $a_1,\dots ,a_n$ are not \alg ically free, we have the set of
cyclic \gra s evaluated at $\a=(a_1,\dots ,a_n)$
\[
(\e_\a)^n(\d{\Cn}) \ .
\]

\subsection{Semicircular systems}

Recall that in a $C^*$-\pr ty space $(A,\var)$ a semicircular system is an
$n$-tuple $(S_1,\dots ,S_n)$ of selfadjoint \ncc random variables
which are freely independent and have (0,1) semicircle distributions.
The von Neumann algebra of $(S_1,\dots ,S_n)$ in the GNS representation
associated with the restriction of $\var$ is isomorphic to the free group
factor $L(F(n))$ and the restriction of $\var$ is the trace-state.

There is a natural realization of a semicircular system on the full Fock space
\[
{\T}({\C}^n) =\bigoplus_{k\geq 0} ({\C}^n)^{\otimes k}
\]
where $ ({\C}^n)^{\otimes 0}={\C}1$ with 1 the vacuum vector.
If $e_j$, $1\leq j\leq n$ are the standard orthonormal basis vectors in
${\C}^n$, let
\[
\ell_j\xi = e_j\otimes\xi \ , \qquad
r_j\xi=\xi\ot e_j
\]
be the left and right creation operators. Then
$s_j=\ell_j+\ell^*_j$ \ $(1\leq j\leq n)$ is a semicircular system w.r.t.~the
vacuum expectation $\<\,\dt 1 ,1\>$ and similarly $d_j=r_j+r^*_j$ \
$(1\leq j\leq n)$.

In particular ${\T}({\C}^n)$ identifies with $L^2(W^*(s_1,\dots ,s_n))$ and
denoting by $P_0,P_1,\dots $ the orthonormal polynomials w.r.t.~a
semicircular distribution we have
\[
e_{i_1}^{\ot k_1}\ot\dots\ot e_{i_p}^{\ot k_p} =
P_{k_1}(s_{i_1})\dots P_{k_p}(s_{i_p})1 =
\ell_{i_1}^{k_1}\dots \ell_{i_p}^{k_p}1
\]
where \ $i_j\neq i_{j+1}$ \ $(1\leq j\leq p-1)$ and $k_j>0$ \
$(1\leq j\leq p)$.

The \po s $P_k(x)$, \ $k\geq 0$ are equal to $C^1_k(x/2)$ where $C^1_k(t)$ are
Gegenbauer \po s, which satisfy the generating function relation
\[
(1-2rt+r^2)^{-1} =\sum_{n\geq 0} C^1_n(t)r^n
\]
(see [9, ch.IX, $\S$4, sec.12]).

\section{Exponentiating noncommutative vector fields}

\bs\noindent{\bf 2.1} \ \ If ${\mathcal E}$ is a ${\Cn}$-bimodule, let
${\VE}$ denote ${\mathcal E}^n$. If $K=(K_j)_{1\leq j\leq n}\in{\VE}$
let $D_K:{\Cn}\to{\mathcal E}$ be the derivation such that
$D_KX_j=K_j$, i.e.,
\[
D_K=\sum_{1\leq j\leq n} m_{K_j}\p_j \ .
\]
In particular if ${\mathcal E}$ is a unital Banach-algebra and the
bimodule structure arises from a unital homomorphism of ${\Cn}$
into ${\mathcal E}$ which takes $X_j$ to $T_j$, then
\[
\frac{d}{d\e} \ P(T_1+\e K_1,\dots ,T_n+\e K_n)|_{\e=0} = D_K P
\]
where \ $P\in {\Cn}$.

\bs\noindent{\bf 2.2} \ \
On ${\Cn}$ we define seminorms
\[
\left| \sum_{p\geq 0} \,\sum_{1\leq j_1,\dots ,j_p\leq n}
c_{j_1,\dots ,j_p} X_{j_1}\dots X_{j_p}\right|_{R,k}
= \sum_{p\geq k} \,\sum_{1\leq j_1,\dots ,j_p\leq n}
|c_{j_1,\dots ,j_p}| R^{p-k}p! ((p-k)!)^{-1} \ .
\]

The completion of ${\Cn}$ w.r.t.~$|\,\dt\,|_{R,j}$, $0\!\leq\! j\!\leq\! k$
will be denoted ${\C}_{\<n\>,R,k}$ or ${\C}\<X_1,\dots ,X_n\>_{\!R,k}$
and identifies with a subalgebra of the \alg \ of \ncc formal power
series\linebreak denoted ${\C}_{\<\<n\>\>}$ or ${\C}\<\<X_1,\dots ,X_n\>\>$.

On \ Vect\,${\C}\<X_1,\dots ,X_n\>$ we have corresponding seminorms
\[
|(K_1,\dots ,K_n)|_{R,k} = \max_{1\leq s\leq n} |K_s|_{R,k} \ .
\] The completions coincide with \ Vect\,${\C}_{\<n\>,R,k}$.

If $K\in {\mb{Vect}}\,{\Cn}$ and $P\in{\Cn}$, then
\[
|D_K P|_{R,0} \leq |K|_{R,0}\, |P|_{R,1} \ .
\]
In particular $D_K P$ can be defined for $K\in  {\mb{Vect}}\,{\C}_{\<n\>,R,0}$
and $P\in {\C}_{\<n\>,R,1}$ as an element in ${\C}_{\<n\>,R,0}$.
Also if $0<R<R'$ then $|P|_{R,1}\leq C|P|_{R',0'}$ for some constant
independent of $P$, so that
$D_KP\in {\C}_{\<n\>,R,0}$ \ if \ $P\in {\C}_{\<n\>,R',0}$ and
$K\in {\mb{Vect}}\,{\C}_{\<n\>,R,0}$. Iterating, we have $D^m_KP\in
{\C}_{\<n\>,R,0}$
\ if
\ $P\in {\C}_{\<n\>,R',0}$ and $K\in {\mb{Vect}}\,{\C}_{\<n\>,R',0}$ for
some
$R'>R$.

\bs\noindent{\bf 2.3} \ \
The map $\Phi_n:{\Cn}\to{\C}_{\<1\>}$ is defined by
\[
\Phi_n\left( \sum_{p\geq 0} \,\sum_{1\leq j_1,\dots ,j_p\leq n}
c_{j_1,\dots ,j_p} X_{j_1}\dots X_{j_p}\right)=
  \sum_{p\geq 0} \,\sum_{1\leq j_1,\dots ,j_p\leq n}
|c_{j_1,\dots ,j_p}| X_1^p \ .
\]
Clearly \ $|\Phi_n(P)|_{R,k} = |P|_{R,k} = (D^k_{X_1}\Phi_n(P))(R)$
where $R>0$, $k\geq 0$.

Obviously $\Phi_n$ extends to a map of ${\C}_{\<n\>,R,k}$ to
${\C}_{\<1\>,R,k}$.

\bs\noindent{\bf 2.4} \ \
If \ $P\in{\Cn}$ let $|P|$ be the \ncc polynomial with coefficients the
absolute values of the coefficients of $P$. Also, if $P,Q\in {\Cn}$ we shall
write $|P|\leq |Q|$ if the inequality holds among the coefficients of $|P|$
and $|Q|$ which are nonnegative numbers).

We extend this definition to \ Vect\,${\Cn}$ by putting
$|(K_j)_{1\leq j\leq n}|=(|K_j|)_{1\leq j\leq n}$ and\linebreak
$|K|\leq |K'|$ if $|K_j|\leq |K'_j|$ $(1\leq j\leq n)$ where
$K=(K_j)_{1\leq j\leq n}$, \ $K'=(K'_j)_{1\leq j\leq n}$.

We shall also write $|P|\vee |Q|$ for the coefficientwise maximum.

\bs\noindent{\bf 2.5} \ \ The analogue of $\Phi_n$ on  Vect\,${\Cn}$
is the map $\Psi_n:{\mb{Vect}}\,{\Cn}\to{\mb{Vect}}\,{\C}_{\<1\>}$
defined by
\[
\Psi_n((K_j)_{1\leq j\leq n}) =\Phi_n(K_1)\vee\dots\vee \Phi_n(K_n) \ .
\]
We have
\[
|K|_{R,k} \leq (D^k_{X_1} \,\Psi_n(K))(R)
\]
with equality if $k=0$.

\bs\noindent{\bf 2.6 Lemma.} \ {\it Let $K,K'\in {\rm{Vect}}\,{\Cn}$, and
$P,P'\in{\Cn}$. Then}
$|D_KP|\leq D_{|K|}|P|$.\\
{\it If $|K|\leq |K'|$ and $|P|\leq |P'|$, then}
\[
D_{|K|}|P|\leq D_{|K'|}|P'| \quad {\mb{and}}
\Phi_n(P)\leq\Phi_n(P') \ , \quad \Psi_n(K)\leq\Psi_n(K') \ .
\]
{\it Moreover we have}
\[
\Phi_n(D_{|K|}|P|) \leq D_{\Psi_n(K)}\Phi_n(P)
\]

The proof reduces to the most obvious majorizations and is left
to the reader.

\bs\noindent{\bf 2.7 Theorem.} \
{\it Let $0<R<R'$ and let $K\in
{\rm{Vect}}\,{\C}_{\<n\>,R',0}$ and $P\in {\C}_{\<n\>,R',0}$.
Then}
\[
|D^m_KP|_{R,0}\leq 1\dt 3\dt 5\dt\dots(2m-1)(R')^{m+1}
(R'-R)^{-2m-1} |K|^m_{R',0} |P|_{R',0} \ .
\]
{\it In particular, $D^m_KP\in {\C}_{\<n\>,R,0}$ and if
$0<r< |K|^{-1}_{R',0}(R'-R)^2(2R')^{-1}$, then}
\[
\sum_{m\geq 0} |D^m_KP|_{R,0} \frac{r^m}{m!} \ < \ \infty
\]
{\it i.e., $P$ is an analytic vector for $D_K$ in} ${\C}_{\<n\>,R,0}$.

\bs
{\bf Proof.}  Clearly all assertions are an easy consequence of the estimate
for $|D^m_KP|_{R,0}$ and it suffices to prove it when $K\in {\mb{Vect}}\,{\Cn}$
and $P\in{\Cn}$.

If \ $|P|_{R',0}=M$ then $(\Phi_n(P))(R')=M$ and the coefficient of
$X^k_1$ of $\Phi_n(P)$ is then majorized by $M(R')^{-k}$ so that
\[
\Phi_n(P)\leq M(1-X_1/R')^{-1} \ .
\]
Similarly, if $|K|_{R',0}=N$ we have
\[
\Psi_n(K)\leq N(1-X_1/R')^{-1} \ .
\]

Using Lemma 2.6 we get
\[
D^m_{\Psi_n(K)}\Phi_n(P)\leq
|K|^m_{R',0} \, |P|_{R',0}\, D^m_{\L}\L
\]
where \ $\L=(1-X_1/R')^{-1}$. By induction we easily get
\[
D^m_{\L}\L = 1\dt 3\dt 5\dt\dots(2m-1)(R')^{-m} \L^{2m+1} \ .
\]
It follows that
\begin{eqnarray*}
|D^m_KP|_{R,0} &\leq & |K|^m_{R',0} \, |P|_{R',0}\, D^m_{\L}\L(R)\\
&=& (2m)!(m!)^{-1}\, 2^{-m}(R')^{m+1} (R'-R)^{-2m-1}
|K|^m_{R',0} \, |P|_{R',0}
\end{eqnarray*}
\qed

\section{Cyclic gradients and exponentiation}

\subsection{The basic property of cyclic gradients}
Let $(M,\tau)$ be a von Neumann \alg \ with a normal faithful trace
state, i.e., a tracial $W^*$-\pr ty space, and let $T_j\in M$,
$K_j\in M$, $1\leq j\leq n$. If $P\in{\Cn}$, then
\[
\frac{d}{d\e}\,\tau(P(T_1+\e K_1,\dots ,T_n+\e K_n))|_{\e =0}=
\sum_{1\leq j\leq n} \tau((\d_j P)(T_1,\dots ,T_n)K_j) \ .
\]
Thus the cyclic gradient $(\d P)(T_1,\dots ,T_n)$ is precisely the \gra \ at
$(T_1,\dots ,T_n)$ of the function
\[
M^n\ni (T_1,\dots ,T_n)\to\tau(P(T_1,\dots ,T_n))\in{\C}
\]
when we use the bilinear scalar product given by
$\sum_{1\leq j\leq n} \tau(a_jb_j)$ on
$M^n={\mb{Vect}}\,M$ (w.r.t.~the bimodule  structure induced by
$\e_{(T_1,\dots ,T_n)}$). Under the sesquilinear scalar product
$\sum_{1\leq j\leq n} \tau(a_jb^*_j)$ the gradient would be
$(\d P(T_1,\dots ,T_n))^*=((\d_j P(T_1,\dots ,T_n))^*)_{1\leq j\leq n}$.

\bs\noindent{\bf 3.2} \ \
On ${\Cn}$ we consider the involution $P\to P^*$ such that $X_j=X^*_j$, i.e.,
$(cX_{i_1}\dots X_{i_p})^*={\bar c}X_{i_p}\dots X_{i_1}$.
If $T_j\in M$ then
\[
(P(T_1,\dots ,T_n))^* = (P^*(T^*_1,\dots ,T^*_n)
\]
and if $T_j=T^*_j$, then
\[
(P(T_1,\dots ,T_n))^* = (P^*(T_1,\dots ,T_n) \ .
\]
It is also easily seen that $\d_jP^*=(\d_jP)^*$. Similarly
$\p_jP^*=\widetilde{(\p_j P)}^*$ where on
${\Cn}\ot{\Cn}$ we define
$(\xi\ot\eta)^*=\xi^*\ot\eta^*$ and $\widetilde{\xi\ot\eta}=\eta\ot\xi$.

\bs\noindent{\bf 3.3} \ \
Let $Y_j=Y^*_j\in M$, $1\leq j\leq n$ be \alg ically free so that ${\Cn}$
and ${\C}\<Y_1,\dots ,Y_n\>$ are isomorphic, via $\e_{(Y_1,\dots ,Y_n)}$.
This turns $L^2(M,\tau)$ into a ${\Cn}$ bimodule. Assume also
$\{Y_1,\dots ,Y_n\}$ generates $M$. If $K\in {\mb{Vect}}\,L^2(M,\tau)$ we
shall consider $D^\e_K=D_K\circ\e^{-1}_{(Y_1,\dots ,Y_n)}$
which is a derivation of  ${\C}\<Y_1,\dots ,Y_n\>$  into $L^2(M,\tau)$.
We shall view $D_K^\e$ as an unbounded densely defined operator
$L^2(M)$ with domain of definition  ${\C}\<Y_1,\dots ,Y_n\>$.

\bs\noindent{\bf 3.4 Proposition.} \
{\it The following are equivalent conditions on $K$.}
\begin{description}
\item{\rm{(i)}} \ $\sum_{1\leq j\leq n}\tau((\d_jP)(Y_1,\dots ,Y_n)K_j)=0$ \
{\it for all} $P\in {\Cn}$.
\item{\rm{(ii)}}  \ $\tau(D^\e_KQ)=0$ {\it for all} $Q\in{\C}\<Y_1,\dots
,Y_n\>$
\item{\rm{(iii)}}  {\it If \ $P_1,P_2\in {\C}\<Y_1,\dots ,Y_n\>$ \ then in
$L^2(M,\tau)$, \ $\<D^\e_KP_1,P_2\>=-\<P_1,D^\e_{K^*}P_2\>$\\
where} $K^*=(K^*_1,\dots ,K^*_n)$.
\end{description}

\bs
{\bf Proof.} (i)$\,\Leftrightarrow\,$(ii) follows from 2.1 and 3.1. We
have
\begin{eqnarray*}
\<D^\e_KP_1,P_2\> + \<P_1,D^\e_{K^*}P_2\> &=&
\tau((D^\e_KP_1)P_2^* + P_1(D^\e_{K^*}P_2)^*) \\
&=&\tau((D^\e_KP_1)P^*_2 +P_1(D^\e_K P^*_2)) \\
&=&\tau(D^\e_K(P_1P^*_2)) \ .
\end{eqnarray*}
Clearly $\tau(D_K(P_1P^*_2))=0$ for all $P_1,P_2\in{\C}\<Y_1,\dots ,Y_n\>$
is equivalent to (ii)\\ so that (ii)$\,\Rightarrow\,$(iii). \qed

\bs\noindent
{\bf 3.4 Corollary.} \ {\it  Assume $K$ satisfies the equivalent conditions
of Proposition $3.3$. Then:}
\begin{quote} {\rm{a)}} \ $D^\e_K$ {\it is closable}\\
{\rm{b)}} \ {\it If $K=K^*$, the densely defined operator $D^\e_K$ is
antisymmetric.}\end{quote}

\bs\noindent{\bf 3.5} \ \ {\it If $K$ satisfies the equivalent conditions
of Proposition $3.3$, we shall say $K$ is a trace-preserving or
$\tau$-preserving (if we want to specify the trace) \ncc vector field.}

Note also that the equivalence in Proposition 3.3 actually holds
more generally for $K\in {\mb{Vect}}\,L^1(M,\tau)$ while (iii) can be
adapted to this case in terms of the duality of $L^1(M,\tau)$ and~$M$.

\bs\noindent{\bf 3.6} \ \
If $K=K^*\in {\mb{Vect}}\,M$ then $D^\e_K P\in M$ if $P\in
{\C}\<Y_1,\dots ,Y_n\>$ 	and $D^\e_KP^*=(D^\e_KP)^*$ so that  $D^\e_K$
is a symmetric derivation (see 3.2.21 in [3]) of
$C^*(Y_1,\dots ,Y_n)$.

If additionally $K$ is in ${\C}\<Y_1,\dots ,Y_n\>$  and is
trace-preserving then by Corollary 3.4,\linebreak $H=-iD^\e_K$ is a symmetric
unbounded operator defined on ${\C}\<Y_1,\dots ,Y_n\>$ which is a dense
subspace in $L^2(M,\tau)$. Since $D^\e_K$ is a derivation, we have
\[
D^\e_KP \ = \ i[H,P]
\]
if \ $P\in {\C}\<Y_1,\dots ,Y_n\>$. This makes $D^\e_K$ a
{\it spatial derivation of} $C^*(Y_1,\dots ,Y_n)$ {\it implemented by} $H$
(see Definition 3.2.54 in [3]).
Summarizing, we have proved

\bs\noindent
{\bf Corollary.} \ {\it  Assume $K=K^*\in {\rm{Vect}}\,{\C}\<Y_1,\dots
,Y_n\>
$ is trace-preserving, then $D^\e_K$ is a spatial derivation of
$C^*(Y_1,\dots ,Y_n)$  implemented on $L^2(M,\tau)$.}

\bs\noindent{\bf 3.7} \ \
Combining the preceding corollary with results on exponentiation of
derivations to automorphisms, one gets results about when $D^\e_K$ can be
exponentiated. For instance, using Proposition 3.2.58 in [3] we get the
following.

\bs\noindent
{\bf Corollary.} \ {\it   Assume $K=K^*\in
{\rm{Vect}}\,{\C}\<Y_1,\dots ,Y_n\> $
is trace-preserving. If  ${\C}\<Y_1,\dots
,Y_n\>\subset L^2(M,\tau)$ consists of analytic vectors for $D^\e_K$, viewed
as an unbounded operator in $L^2(M,\tau)$, then   $D^\e_K$
exponentiates to a one-parameter automorphism group $\exp(tD^\e_K)$ of $M$
which is trace-preserving} $\tau\circ\exp(tD_K)=\tau$.

\bs\noindent{\bf 3.8} \ \ Let $R\geq\|Y_j\|$, $1\leq j\leq n$.
Then $|P(Y_1,\dots ,Y_n)|_2\leq \|P(Y_1,\dots ,Y_n)\|\leq
|P(X_1,\dots ,X_n)|_{R,0}$. Then Theorem 2.7 implies that for some $r>0$
we have
\[
\sum_{m\geq 0} |(D^\e_K)^m P|_2\,\frac{r^m}{m!} \ < \ \infty
\]
if $P\in{\C}\<Y_1,\dots ,Y_n\>$, i.e., $P$ is an analytic vector in $L^2(M)$.
Thus we have actually proved:

\bs
{\bf  Theorem.} \ {\it
  Assume $Y_j=Y^*_j\in M$,
$1\leq j\leq n$, and $K=K^*\in{\rm{Vect}}\,{\C}\<Y_1,\dots ,Y_n\>$.\linebreak
Assume $\{Y_1,\dots ,Y_n\}$ generates $M$ and the $Y_j$'s are \alg ically
free. If $K$ is trace-preserving, then $D^\e_K$ exponentiates to a
one-parameter group of trace-preserving automorphism of $M$.}

\bs\noindent{\bf 3.9} \ \ Let $R\geq\|Y_j\|$, $1\leq j\leq n$.
We could have replaced during all steps the ring
${\C}\<Y_1,\dots ,Y_n\>$ with the ring
\[
\e_{(Y_1,\dots ,Y_n)}{\C}_{\<n\>, > R}
\]
of power-series of radius of convergence $>R$. The only additional assumption
is that \alg ic freeness has to be replaced with the stronger requirement
that
\[
\e_{(Y_1,\dots ,Y_n)}{\C}_{\<n\>, > R} \ \to \ M
\]
is one-to-one. With these amendments trace-preserving self-adjoint
$K\in {\mb{Vect}}\, \e_{(Y_1,\dots ,Y_n)}{\C}_{\<n\>, > R}$ also
``exponentiate" to one-parameter groups of automorphisms of $M$.

Of course exponentiate both here and in 3.7, 3.8 means either that we take a
weak closure of the exponential or use the automorphism group implemented by
the unitary group $\exp(it{\overline H})$ where ${\overline H}$ is
the closure of $H$.

\section{Endomorphic orbits}

\bs\noindent{\bf 4.1} \ \ Let $Y_j=Y^*_j\in (M,\tau)$, $1\leq j\leq n$.
We shall denote by $M_h=\{a\in M\mid a=a^*\}$ the hermitian elements of $M$
and by $L^p(M_h,\tau)$ the corresponding $L^p$-space. The scalar product
$\sum_{1\leq j\leq n}\tau(a_jb^*_j)$ on $(L^2(M,\tau))^n$ restricts to a real
scalar product on  $(L^2(M_h,\tau))^n$. Note also that the set of cyclic \gra s
is self-adjoint when evaluated at $Y_1,\dots ,Y_n$. The {\it orthogonal to
cyclic \gra s at} $Y_1,\dots ,Y_n$ is:
\[
(\d {\C}^{\perp}_{\<n\>})(Y_1,\dots ,Y_n) =\{
Z\in (L^2(M_h,\tau))^n\mid \<Z,\d P(Y_1,\dots ,Y_n)\>=0 \ \ {\mb{if}} \ \
P\in{\Cn}\} \ .
\]
Clearly in the previous formula we might consider only $P=P^*$. Also, we have
corresponding $p$-spaces replacing $L^2$ by $L^p$.

\bs\noindent{\bf 4.2} \ \ By Aut$(M,\tau)$ (resp.~End$(M,\tau)$)
we denote auto\mor s (resp.~unital endo\mor s) $\a$ of $M$ such that
$\tau\circ\a=\tau$. If $M$ is a II${}_1$-factor then preservation of the trace
is automatic and we can write Aut$(M)$  (resp.~End$(M)$).
{\it If \ $Y\in M^n_h$ \ let} $AO(Y)={\mb{Aut}}(M,\tau)\dt Y$
{\it (resp.}~$EO(Y)={\mb{End}}(M,\tau)\dt Y$) {\it be the automorphic
(resp.~endomorphic) orbit of} $Y$.

We define the tangent set to $AO(Y)$ (resp.~$EO(Y)$) at $Y$ to be
\begin{eqnarray*}
T\!AO(Y) &=& \{\xi\in (L^2(M_h,\tau))^n\mid\exists \ \eta_k\in AO(Y), \\
&&\quad \exists \ \mu_k\in{\R}\backslash\{0\}, \
\mu_k\to 0, \ \lim_{k\to\infty}|\mu^{-1}_k(\eta_k-Y)-\xi|_2=0\}
\end{eqnarray*}
(resp.~$TEO(Y)$ defined similarly using EO instead of AO).

There are several variants of these sets $T\!AO_p(Y)$,
$wT\!AO_p(Y)$ (resp.~$TEO_p(Y)$, $wTEO_p(Y)_n$) where we take
$x\in (L^p(M_h,\tau))^n$ and require the limit to hold in $p$-norm
$(1\leq p\leq\infty)$ or for $wT\!AO_p(Y)$ (resp.~$TEO_p(Y)$) to hold
w.r.t.~the
weak convergence in the duality with $(L^q(M_h,\tau))^n$
$(p^{-1}+q^{-1}=1)$. If $p=2$ we simply write $wT\!AO(Y)$ (resp.~$wTEO(Y)$).
Remark also that $wTEO(Y)\supset TEO(Y)\cup T\!AO(Y)\cup wT\!AO(Y)$.

\bs\noindent
{\bf 4.3 Proposition.} \
  \ $wTEO(Y)\subset (\d{\C}^{\perp}_{\<n\>})(Y)$.

\bs
{\bf Proof.} Let \ $\eta_k\in EO(Y)$ so that $w-\lim_{k\to\infty}
(\mu^{-1}_k(\eta_k-Y)-\xi)=0$ in $(L^2(M_h,\tau))^n$. Since the weak
convergence is for a sequence, we have uniform boundedness in $L^2$ and hence
$\mu^{-1}_k=O(|\eta_k-Y|^{-1}_2)$.  On the other hand, since $\eta_k,Y$
are in $(M_h)^n$, we have
\[
\tau(P(\eta_k)-P(Y)-\sum_j(\eta_k-Y)_j(\d_j P)(Y)) =O(|\eta_k-Y|^2_2)
\]
where \ $P\in{\Cn}$.  Since
\begin{eqnarray*}
\mu^{-1}_k O(|\eta_k-Y|^2_2)& =& O(|\eta_k-Y|^{-1}_2 |\eta_k-Y|^2_2)\\
  &=& O(|\eta_k-Y|_2) = o(1)
\end{eqnarray*}
and $\tau(P(\eta_k))=\tau(P(Y))$ we have
\[
\lim_{k\to\i}\tau\left(\sum_j\mu^{-1}_k(\eta_k-Y)_j
(\d_jP)(Y)\right) = 0 \ .
\]
We infer \ $\sum_j\tau(\xi_j(\d_jP)(Y))=0$, i.e.,
$\xi\in(\d{\C}^{\perp}_{\<n\>})(Y)$. \qed

\bs\noindent{\bf 4.4} \ \ If $Y\in (M_h)^n$ is a generator of $M$,
then $EO(Y)$ parametrizes End$(M,\tau)$.

{\it It is also easily seen that $EO(Y)$ is a closed set in the topology of
the $L^2$-metric.}

On the other hand, if $Y=(Y_1,\dots ,Y_n)$ {\it is a generating
semicircular system of $L(F(n))$, then $AO(Y)$ being $L^2$-closed
would imply the non-iso\mor \ of free group factors.}

Indeed $L(F(\i))$ has a sequence of auto\mor s $\a_k$ so that
$\a_k(\l(g_p))=\l (g_q)$ if $q-p\equiv 1$ (mod $k$) when $1\leq p\leq k$
and $p=q$ if $p>k$, which converge pointwise on $L^2$ to a non-invertible
endo\mor . Thus the iso\mor \ of $L(F(n))$ and $L(F(\i))$ would imply the
automorphic orbit of $Y$ is not closed.

\bs\noindent{\bf 4.5 Remark.} \ In connection with free \pr ty, there
is a third
kind of orbits which occur: the orbits of the relation of
equivalence in distribution. The {\it distribution orbit of} $Y$
denoted $DO(Y)$ is the set
$\{Y'\in M^n_h\mid Y'$ and $Y$ equivalent in distribution$\}$
where equivalence in distribution means $\tau(Y_{i_1}\dots Y_{i_p})=\tau
(Y'_{i_1}\dots Y'_p)$ for all $p\in {\mathbb N}$ and $i_j\in\{1,\dots ,n\}$.
Equivalently, $Y'\in DO(Y)$ if there is a unital trace-preserving
iso\mor \ $\a:W^*(I,Y_1,\dots ,Y_n)\to W^*(I,Y'_1,\dots ,Y'_n)$ so that
$\a(Y_j)=Y'_j$. Clearly
$DO(Y)\supset EO(Y)$. Then the tangent sets $TDO(Y)$,
$T_pDO(Y)$, $wTDO(Y)$, $wTDO_p(Y)$ are defined similarly to the
endomorphic and automorphic cases and we have the analogue of
Proposition 4.3, i.e.:
\[
wTDO(Y)\subset (\d{\C}^\perp_{\<n\>})(Y) \ .
\]

\section{Real and complex cyclomorphic maps}

\medskip\noindent{\bf 5.1} \ \ Throughout section 5, by $(M,\tau)$,
$(N,\nu)$, $(A,\a)$, $(B,\b)$, $(C,\g)$ we shall denote von Neumann algebras
with normal faithful trace-states.

\bs{\bf Definition.} A differentiable map $f:\O\to{\R}$, where
$\O\subset M^n_h$ is an open set in the norm-topology, is ${\R}$-cyclomorphic
if for every $Y\in\O$ there is $T\in (L^1(M_h,\tau))^n$ such that for $K\in
M^n_h$ we have a
\[
df[Y](K_1,\dots,K_n) =\sum_{1\leq j\leq n} \tau(T_jK_j)
\]
and $T$ is in the $L^1$-norm closure of $\d{\Cn}(Y)\cap M^n_h$.
We shall abbreviate saying $f$ is ${\R}cm$.

\bs{\bf 5.2 \ Definition.} A differentiable map $f:\O\to N^p_h$, where
$\O\subset M^n_h$ is  open in the norm-topology, is ${\R}$-cyclomorphic
if for every ${\R}$-cyclomorphic map $f:\o\to{\R}$, where $\o\subset N^p_h$
is open in the norm-topology, we have that $f\circ F:\O\cap F^{-1}(\o)\to{\R}$
is ${\R}$-cyclomorphic. Also in this case this will be abbreviated by $F$ is
${\R}c m$.

\bs\noindent{\bf 5.3} \ \ Since
${\R}={\C}_h$, a ${\R}c m$ map $f:\O\to{\R}$ is also a map
$f:\O\to{\C}_h$. It is easily seen that definition 5.1 for $f:\O\to{\R}$
and definition 5.2 for $f:\O\to{\C}_h$ are equivalent. For this it suffices to
remark that the ${\R}c m$ maps $g:\o\to{\R}$, $\o\subset{\C}_h$ are
just the differentiable maps and that if $g:\o\to{\R}$ is differentiable,
$f(Y)\in \o$, $Y\in\O$, then $d(g\circ f)[Y]=g'(f(Y))dg[Y]$.

\bs\noindent{\bf 5.4 \ Examples.} \ \ a) If $P=P^*\in{\Cn}$ then the map
$f_P:M^n_h\to{\R}$, defined by $f_P(Y)=\tau(P(Y))$ is ${\R}cm$.
Indeed $df_P[Y](K)=\tau(\sum_j(\d_j P)(Y)K_j)$.

b) If \ $f_1,f_2:\O\to{\R}$, where $\O\subset M^n_h$, are ${\R}c m$,
then also $f_1+f_2$ and $f_1f_2$ are ${\R}c m$. The ${\R}c
m$-maps  from $\O$ to ${\R}$ form an \alg \ over ${\R}$ with unit.

\bs\noindent
{\bf 5.5 \ Lemma.} \ {\it
  A differentiable map $F:\O\to N^p_h$, where
$\O\subset M^n_h$ is ${\R}c m$ iff for every $P=P^*\in{\C}_{\<p\>}$
the map $f_P\circ F$ is ${\R}c m$, where} $f_P: N^p_h\to{\R}$,
$f_P(T)=\nu(P(T))$.

\bs
{\bf Proof.} The ``only if\," is obvious by 5.4a). To prove the ``if\,"-part,
remark that Definition 5.1 for a differentiable map $f:\o\to{\R}$, $\o\subset
N^p$ open, is equivalent to saying that $f$ is ${\R}c m$ iff for every
$T\in N^p_h$ there is a sequence of polynomials $P_s=P^*_s\in{\C}_{\<p\>}$
such that $\|df_{P_s}[T]-df[T]\|\to 0$ as $s\to\i$ in ${\mathcal
L}(N^p_h,{\R})$. If $Y\in\O$, $F(Y)\in\o$, then
$d(f\circ F)[Y]=df[F(Y)]\circ dF[Y]$, so that
\[
\|d(f\circ F)[Y]-d(f_{P_s}\circ F)[Y]\| \leq
\|df[F(Y)]-df_{P_s}[F(Y)]\| \ \|dF[Y]\| \ .
\]
Hence $d(f\circ F)[Y]$ is in the closure of differentials
$d(f_{P_s}\circ F)[Y]$, i.e. in the closure of linear maps given by
cyclic \gra s. Hence $f\circ F$ is  ${\R}c m$. \qed

\bs\noindent{\bf 5.6 \ Examples.} \ \ a) If $P_j=P^*_j\in{\Cn}$,
$1\leq j\leq p$, then the map $F_p:M^n_h\to M^p_h$ given by
$F_P(Y)=(P_j(Y))_{1\leq j\leq p}$ is ${\R}c m$. Indeed if
$Q=Q^*\in {\C}_{\<p\>}$, then
\[
(f_Q\circ F_P)(Y) =\tau(Q(P(Y))=f_{Q\circ P}
\]
and the assertion follows from Lemma 5.5 and Example 5.4a).

b) Let $\O_1\subset A^n_h$, $\O_2\subset B^p_h$, be open sets in the
norm-topology, and let $F:\O_1\to B^p_h$, $G:\O_2\to C^q_h$ be
${\R}cm$ maps, so that $F(\O_1)\subset\O_2$. Then $F\circ G$
is also an ${\R}cm$ map.

c) If \ $f:\O\to{\R}$ and $F:\O\to N^p_h$, where $\O\subset M^n_h$,
are ${\R}cm$, then $H:\O\to N^p_h$ defined as
$H(Y)=f(Y)F(Y)$ is ${\R}cm$. Indeed by Lemma 5.5 it suffices to show
that if $P=P^*\in {\C}_{\<p\>}$ we have $f_P\circ H$ is ${\R}cm$.
Clearly there is no loss of generality to assume $P$ is homogeneous of degree
$q$. Then $(f_P\circ H)(Y)=(f(Y))^q(f_P\circ H)(Y)$ is ${\R}cm$
by 5.4b).

\bs\noindent
{\bf 5.7 \ Proposition.} \ {\it
Let $f:\O\to{\R}$, $\O\subset M^n_h$ be a
${\R}cm$ map and let $(\a_t)_{t\geq 0}\stt {\rm{End}}(M,\tau)$ be a
semigroup, i.e., $\a_0={\rm{id}}_M$, $\a_{s+t}=\a_s\circ\a_t$. Let further
$Y\in\o$ be such that for $t\in [0,1]$, $\a_t(Y)\in\O$,
$t\to\a_t(Y)$ is norm-continuous and the right derivative}
\[
\lim_{t\downarrow 0} t^{-1}(\a_t(Y)-Y) = K\in M^h_n
\]
{\it exists in norm-convergence. Then}
\[
f(\a_t(Y))=f(Y)
\]
{\it for all} $t\in [0,1]$.

\bs
{\bf Proof.} By Proposition 4.3 we have  that $K\in\d{\C}^{\perp}_{\<n\>}(Y)$
which implies $df[Y](K)=0$, i.e. $\lim_{t\downarrow 0} t^{-1}(f(\a_t(Y))
-f(Y)) =0$.
By the semigroup property we have
$\lim_{t\downarrow 0} t^{-1}(\a_{s+t}(Y)-\a_s(Y)) =\a_s(K)$ and the same
argument applied to $\a_s(Y)$, $0\leq s <1$ gives
\[
\lim_{t\downarrow 0} t^{-1}(f(\a_{t+s}(Y)) -f(\a_s(Y)))=0 \  .
\]
Therefore $f(\a_t(Y))$ as a function of $t\in [0,1]$ is continuous and has
zero right derivatives at each point of [0,1). We infer $f(\a_t(Y))$
is constant on [0,1]. \qed

\bs\noindent
{\bf 5.8 Corollary.} \ {\it  \ If \ $f:M^n_h\to{\R}$ is ${\R}cm$ then $f$ is
unitary-invariant, i.e., if $U\in M$ is unitary and $Y=
(Y_j)_{1\leq j\leq n}\in M^n_h$ then}
$f((UY_jU^*)_{1\leq j\leq n})=f(Y)$.

\bs\noindent
{\bf 5.9 \ Proposition.} \ {\it  Let $\O\stt M_h$ be open in the norm-topology.
A differentiable map\linebreak $f:\O\to{\R}$ is ${\R}cm$ iff
$df[Y](m)=\tau(\xi m)$ for some $\xi=\xi^*\in L^1(W^*(Y),\tau)$
for each} $Y\in\O$.

\bs
{\bf Proof.} This follows from the fact that
$\d{\C}_{\<1\>}(Y)\cap M_h= M_h\cap{\C}_{\<1\>}(Y)$ is a dense
subset of $L^1(W^*(Y)_h,\tau)$, which is clearly a closed subset of
${\mathcal L}(M_h,{\R})$ via the identification $\xi\rightsquigarrow
\tau(\xi\dt)$. \qed

\bs\noindent
{\bf 5.10 \ Corollary.} \ {\it  Let} $M_h^{\mbox{\rm\footnotesize{inv}}}$
{\it be the elements with bounded inverse in $M_h$. The map}
\[
M_h^{\mbox{\rm\footnotesize{inv}}}\ni Y \ \to \ Y^{-1}\in M_h
\]
{\it is} ${\R}cm$.

\bs
{\bf Proof.} By 5.5 it suffices to check that
\[
M_h^{\mbox{\rm\footnotesize{inv}}}\ni Y \ \to \ \tau( Y^{-k})\in {\R}
\]
is ${\R}cm$ for all $k>0$. Since the differential is
\[
M_h\ni K \ \to \ -k\tau(Y^{-k-1}K)
\]
the assertion follows from Proposition 5.9.

\bs\noindent
{\bf 5.11 \ Proposition.} \ {\it  Let $F:\O\to N^p_h$ be a
differential map, where $\O\stt M^n_h$ is open in the norm-topology.
If for each $Y$ there are ${\R}cm$ maps $F_k:\O\to N^p_h$ such that
for} $k\to\i$,
\begin{eqnarray*}
&&\|F_k(Y)-F(Y)\| \ \to \ 0\\
&&\|dF_k[Y]-dF[Y]\| \ \to \ 0
\end{eqnarray*}
{\it then $F$ is} ${\R}cm$.

\bs
{\bf Proof.} We shall use Lemma 5.5. Let $P\in{\C}_{\<p\>}$ and
$f_P:N^p\to {\R}$, $f_P(Y)=\nu(P(Y))$. We must prove that
$f_P\circ F$ is ${\R}cm$. Remark that for $k\to\i$ we have
\begin{eqnarray*}
&&\|(f_P\circ F)(Y)-(f_P\circ F_k)(Y)\| \ \to \ 0\\
&&\|d(f_P\circ F)[Y]-d(f_P\circ F_k)[Y]\| \ \to \ 0
\end{eqnarray*}
where for the second equality we used the fact that
\[
\|df_P[F(Y)]-df_P[F(Y_k)]\| \ \to \ 0 \ .
\]
Since $d(f_P\circ F_k)[Y]$ converges to $ d(f_P\circ F)[Y]$ we infer
$f_P\circ F$ is  ${\R}cm$. \qed

\bs\noindent
{\bf 5.12 \ Corollary.} \ {\it   Let $\O=\{Y\in M^n_h\mid\|Y_j\|<R, \ 1\leq
j\leq n\}$ and let $G_j=G^*_j\in {\C}_{\<n\>,R,0}$ be a \ncc power series.
Then the map $F_G$ defined by $\O\ni Y\to (G_j(Y))_{1\leq j\leq p}\in
M^p_h$ is} ${\R}cm$.

\bs\noindent
{\bf 5.13 \ Theorem.} \ {\it  Let $A\stt B$ be a unital inclusion of
von Neumann
\alg s and let $Z=(Z_1,\dots ,Z_m)\stt B^m_h$. Let $\Phi_Z:A^n_h\to B_h^{n+m}$
be the map defined by $\Phi_Z(Y_1,\dots ,Y_n)=(Y_1,\dots ,Y_n,
Z_1,\dots ,Z_m)\in B_h^{n+m}$. If $\{Z_1,\dots ,Z_m\}$ and $A$ are independent
or freely independent in $(B,\b)$, then $\Phi_Z$ is} ${\R}cm$.

\bs
{\bf Proof.} Using Lemma 5.5 it suffices to prove under the given assumptions
that $f_P\circ \Phi_Z$ is ${\R}cm$ where $P=P^*\in{\C}_{\<n+m\>}$.
In both cases $\tau(P(Y_1,\dots ,Y_n,Z_1,\dots ,Z_m)$ is a \po \ in the
moments of $(Y_1,\dots ,Y_n)$ and the moments of $(Z_1,\dots ,Z_m)$.
The moments of $(Z_1,\dots ,Z_m)$ are constants, so $f_P\circ\Phi_Z$ is a \po
\ in the moments of $(Y_1,\dots ,Y_m)$, or being a real \po \ in the real and
imaginary parts of such moments it is ${\R}cm$ by 5.4b). \qed

\bs\noindent
{\bf 5.14 \ Proposition.} \ {\it  Let $Y_j\in M_h$, $Z_j\in M_h$, $1\leq j\leq
n$ be such that
$\{Y_1,\dots ,Y_n\}$ and
$\{Z_1,\dots ,Z_n\}$ are independent or freely independent. Let $W_{Y+Z}$
and $W_Z$ be the von Neumann sub\alg s of $M$ generated by
$\{I,Y_1+Z_1,\dots Y_n+Z_n\}$ and respectively $\{I,Z_1,\dots ,Z_n\}$
and let $E$ be the conditional expectation of $M$ onto $W_{Y+Z}$. Then, if
$\xi\in ((W_Z,\tau))^n_h$ is such that $\tau(\sum_{1\leq j\leq n}\xi_j\d_j
P(Z))=0$ for all $P\in{\Cn}$, then $\tau(\sum_{1\leq j\leq n}\xi_j\d_j
P(Y+Z))=0$ for all  $P\in{\Cn}$. In particular
$(E\xi_j)_{1\leq j\leq n}$ is in the orthogonal to}
$\d{\Cn}(Y+Z)$ in $(W_{Y+Z,h})^n$.

\bs
{\bf Proof.} The map $S:M^{2n}_h\to M^n_h$, \
$S((R_j)_{1\leq j\leq 2n}=(R_j+R_{n+j})_{1\leq j\leq n}$ is ${\R}cm$
and hence $S\circ \Phi_Z=\Psi_Z$ is ${\R}cm$, where $\Phi_Z$ is the map in
Theorem 5.13. We have
\[
\Psi_Z:M^n_h\to M^n_h \ , \quad
\Psi_Z(Y_1,\dots ,Y_n)=(Y_j+Z_j)_{1\leq j\leq n} \ .
\]
If $P=P^*\in {\Cn}$ then $d(f_P\circ\Psi_Z)[Y]$ is of the form
$\sum_j\tau(T_j\dt)$ with $(T_j)_{1\leq j\leq n}$ in the $L^1$ closure of
$\d{\Cn}(Y+Z)$ (the $L^1$ closure is actually superfluous here).
Since $\xi\in {\mb{Ker}}\, d(f_P\circ\Psi_Z)[Y]$ we  infer
$d\Psi_Z[Y](\xi)\in {\mb{Ker}}\, df_P[Y+Z]$ for all $P=P^*\in{\Cn}$.
Since $d\Psi_Z[Y](\xi)=\xi$ this yields the assertion of the proposition.
\qed

\bs\noindent{\bf 5.15} \ \ On $M^n_h$ we define the ${\R}c
m$-equivalence relation, by $Y\sim Y'$, if $f(Y)=f(Y')$ for all ${\R}cm$
maps
$f:M^n_h\to{\R}$. In particular taking $f=f_P$ where
$P=P^*\in{\Cn}$ we see that $Y\sim Y'$ implies that
$(Y_1,\dots ,Y_n)$ and $(Y'_1,\dots ,Y'_n)$ are equivalent in distribution.
This means that $Y\sim Y'$ implies there is an iso\mor \
$\rho:W^*(I,Y_1,\dots ,Y_n)\to W^*(I,Y'_1,\dots ,Y'_n)$ such that
$\rho(Y_j)=Y'_j$ and which is trace-preserving. On the other hand,
5.7 and 5.8 give some partial converses.

\bs\noindent{\bf 5.16} \ \ A $k+1$-times differentiable map $F:\O\to N^p_h$,
where $\O\stt M^n_h$ is open in the norm-topology, is $k-{\R}cm$
$(k\in {\mathbb N}\cup\i)$ if for all $i\leq j\leq k$ the $j$-th order
differential $d^jF:\O\x M^j_h\to N^p_h$ is ${\R}cm$.

\bs\noindent{\bf 5.17 \ Examples.} \ \ The maps $F_P$ defined in 5.6a)
are $\i-{\R}cm$. Indeed $d^jF_P:M_h^{j+1}\to M^p_h$ are also maps
of the type defined in 5.6a).

\bs\noindent{\bf 5.18} \ \ We shall now take a look at some weakenings of the
requirement of cyclomorphy.

\medskip
{\bf Definition.} Let $\O\stt M^n_h$ be open. A map
$F:\O\to N^p_h$ is weakly ${\R}cm$ (abbreviated $w{\R}cm$)
at a given point $Y\in\O$ if $F$ is differentiable at $Y$ and
\[
(dF[Y])(\d{\C}^{\perp}_{\<n\>}(Y)\cap M^n_h)\stt
\d{\C}^{\perp}_{\<n\>}(F(Y))\cap N^p_h \ .
\]
In particular a map $f:\O\to{\R}$ is $w{\R}cm$ if $f$ is
  differentiable at $Y$ and
\[
\d{\C}^{\perp}_{\<n\>}(Y)\cap M^n_h\stt {\mb{Ker}} \ df \ .
\]
A map $f:\O\to {\R}$ is pseudo ${\R}cm$ (abbreviated $\psi{\R}cm$)
at $Y\in\O$ if $f$ is continuous at $Y$ and $f(Y+T)=f(Y)+o(\|T\|)$
for $T\in(Y+\d{\C}^{\perp}_{\<n\>}(Y)\cap M^n_h)\cap\O$.

\bs\noindent{\bf 5.19} \ \ The difference between global
$w{\R}cm$ and ${\R}cm$ has to do with certain ultraweak continuity
requirements for the differential. We record these facts as the next
proposition, the easy proof of which is left to the reader.

\bs
{\bf Proposition.} \ {\it  Let $\O\stt M^n_h$ be open and let
$F:\O\to N^p_h$ and $f:\O\to{\R}$ be differentiable maps. Then $f$ is ${\R}cm$
iff $f$ is $w{\R}cm$ at $Y$ and $df[Y]$ is ultraweakly continuous for each
$Y\in\O$. If $F$ is ${\R}cm$ then $F$ is $w{\R}cm$ at each point
$y\in\O$. If $F$ is $w{\R}cm$ at $Y$ and $dF[Y]$ is ultraweakly continuous
at each point $Y\in\O$, then $F$ is} ${\R}cm$.

\bs\noindent{\bf 5.20} \ \ We pass now to our brief discussion of complex
cyclomorphy. We shall need to consider spaces of cyclic \gra s evaluated at
non-selfadjoint $n$-tuples $Z=(Z_j)_{1\leq j\leq n}\in M^n$.

\bs{\bf Definition.} Let $\Theta\stt M^n$ be an open set.
A complex-analytic map $h:\Theta\to{\C}$ is complex cyclomorphic
(abbreviated ${\C}cm$) if for each $Z\in\O$ there is
$S\in (L^1(M,\tau))^n$ such that for $W\in M^n$ we have
$dh[Z](W_1,\dots ,W_n)=\sum_{1\leq j\leq n}\tau(S_jW_j)$ and $S$ is in the
$L^1$-norm-closure of $\d{\Cn}(Z)\cap M^n$. A complex-analytic map
$H:\Theta\to N^p$ is  ${\C}cm$ if for every $h:\theta\to{\C}$,
where $\theta\stt N^p$ is open, we have that $h\circ H:\Theta\cap
h^{-1}(\theta)\to{\C}$ is   ${\C}cm$.

\bs\noindent{\bf 5.21} \ \ Many of the results for ${\R}$-cyclomorphic
maps have correspondents for ${\C}$-cyclomorphic maps and the proofs are
usually analogous. We shall briefly state, without proofs, a few of these.

\bs{\bf Proposition.}  \ {\it A complex analytic map $H:\Theta\to N^p$, where
$\Theta\stt M^n$ is open, is ${\C}cm$ iff for every $P\in
{\C}_{\<p\>}$ the map $h_p\circ F$ is  ${\C}cm$, where}
$h_P:N^p\to{\C}$, $h_P(T)=\nu(P(T))$.

\bs\noindent{\bf 5.22 Examples.} \ a) If $P\in{\Cn}$, then
$h_P:M^n\to{\C}$ defined by $h_P(Z=\tau(P(Z))$ is ${\C}cm$.

b) If $h\in\O\to{\C}$ is ${\C}cm$, $H:\O\to N^p$ is
${\C}cm$, where $\O$ is open, then also $Hh:\O\to N^p$
is  ${\C}cm$. In particular the ${\C}cm$ maps from $\O$ to
${\C}$ form an \alg \ over ${\C}$ with unit.

c) If $P_j\in {\Cn}$, $1\leq j\leq p$, then the map
$H_P:M^n\to M^p$, given by $H_P(Z)=(P_j(Z))_{1\leq j\leq p}$ is
${\C}cm$.

d) Let $\Theta_1\stt A^n$, $\Theta\stt B^p$ be open sets, and let
$H:\Theta_1\to B^p$, $K:\Theta_2\to C^q$ be ${\C}cm$ maps such that
$H(\Theta_1)\stt \Theta_2$. Then $K\circ H$ is  ${\C}cm$.

\bs\noindent {\bf 5.23 \ Proposition.}  \ {\it Let $\Theta\stt M$ be open.
A complex analytic map $h:\Theta\to{\C}$ is ${\C}cm$ iff
$dh[Z](T)=\tau(\xi T)$ for some $\xi\in L^1$-closure of}
${\C}_{\<1\>}(Z)$.

\bs\noindent{\bf 5.24 \ Proposition.} \  {\it Let $M^{\rm{\footnotesize{inv}}}$
be the elements with bounded inverse in $M$. The map}\\
$M^{\rm{\footnotesize{inv}}}\ni Z\to Z^{-1}\in M$
is ${\C}cm$.

\bs\noindent{\bf 5.25 \ Proposition.} \ {\it Let $H:\Theta\to N^p$ be a
complex-analytic map, where $\Theta\stt M^n$ is open. If for each $Z\in M^n$
there are
${\C}cm$ maps $H_k:\Theta\to N^p$ such that for} $k\to\i$
\begin{eqnarray*}
&&\|H_k(Z)-H(Z)\|\to 0 \\
&&\|dH_k[Z]-dH[Z]\| \to 0
\end{eqnarray*}
{\it then $H$ is} ${\C}cm$.

\bs\noindent{\bf 5.26 \ Proposition.} \ {\it Let $\Theta=\{Z\in
M^n\mid\|Z_j\|<R,
\ 1\leq j\leq n\}$ and let $G_j\in {\C}_{\<n\>,R,0}$ be \ncc power series.
Then the map $H_G$ defined by}
\[
\Theta\ni Z\to (G_j)(Z))_{1\leq j\leq p}\in M^p
\]
  {\it is} ${\C}cm$.

\bs\noindent{\bf 5.27 \ Theorem.} \ {\it
Let $A\stt B$ be a unital inclusion of von Neumann \alg s and let
$W=(W_1,\dots ,W_m)\in B^m$. Let $\Xi_W:A^n\to B^{n+m}$ be the
map defined by  $\Xi_W(Z_1,\dots ,Z_n)=(Z_1,\dots ,Z_n,W_1,\dots
,W_m)\in  B^{n+m}$. If $\{ W_1,\dots ,W_m\}$ and $A$ are independent or
freely independent in $(B,\b)$ then $\Xi_W$ is} ${\C}cm$.

\bs\noindent
{\bf 5.28 \ Definition.} A complex-analytic map $H:\Theta\to N^p$, where
$\Theta\stt M^n$ is open, is $k-{\C}cm$ $(k\in{\mathbb N}\cup\i)$
if for all $1\leq j\leq k$, the $j$-th order complex differential
$d^jH:\O\x M^j\to N^p$ is ${\C}cm$.

\bs\noindent{\bf 5.29 \ Example.} If $P_i\in{\Cn}$, $1\leq i\leq p$,
then the map $H_P:M^n\to M^p$ given by $H_P(Z)=(P_i(Z))_{1\leq i\leq p}$
is $\i-{\C}cm$.

\bs\noindent{\bf 5.30} \ \  To define weak ${\C}cm$ maps we first
must adapt some  notation. For non-selfadjoint $Z\in M^n$ we shall denote
by
$\d{\C}^{\perp}_{\<n\>}(Z)$ the set
\[
\{T\in (L^1(M,\tau))^n\mid\sum_{1\leq j\leq n}
\tau(T_j(\d_j P)(Z))=0 \ {\mb{for all}} \ P\in{\Cn}\} \ .
\]
In case there may be some confusion with the notation in the selfadjoint
case, we may emphasize the non-selfadjointness by writing
$(\d{\C}^{\perp}_{\<n\>})_{nh}(Z)$.

\bs{\bf Definition.} Let $\Theta\subset M^n$ be open.
A map $H:\Theta\to N^p$ is weakly ${\C}cm$ (abbreviated
$w{\C}cm$) at a point $Z\in\Theta$ if $H$ is ${\C}$-differentiable at
$Z$ and its ${\C}$-differential satisfies
\[
(dH[Z])(\d{\C}^{\perp}_{\<n\>}(Z)\cap M^n)\stt
\d{\C}^{\perp}_{\<p\>}(H(Z))\cap N^p \ .
\]
In particular a map $h:\Theta\to{\C}$ is $w{\C}cm$ at $Z$ if
$h$ is  ${\C}$-differentiable at
$Z$ and $\d{\C}^{\perp}_{\<n\>}(Z)\cap M^n\stt {\mb{Ker}}\,df$.

\bs\noindent{\bf 5.31 \ Proposition.} {\it Let $\Theta\stt M^n$ be open and let
$H:\Theta\to N^p$ and $h:\Theta\to{\C}$ be maps. Then $h$ is
${\C}cm$ iff $h$ is $w{\C}cm$ at $Z$ and $dh[Z]$ is
ultraweakly continuous for each $Z\in\Theta$. If $H$ is ${\C}cm$
then $H$ is $w{\C}cm$ at each point $Z\in\Theta$. If $H$
  is $w{\C}cm$ at $Z$ and $dH[Z]$ is ultraweakly continuous
at each point $Z\in\Theta$, then $H$ is} ${\C}cm$.

\bs\noindent{\bf 5.32 \ ${\R}cm$ and ${\C}cm$ manifolds.}
If $\O_1$ and $\O_2$ are open in $M_h^n$ (resp.~in $M^n$) then a
${\R}cm$-iso\mor \ (resp.~${\C}cm$ iso\mor) is a
${\R}cm$ (resp.~${\C}cm$) map $F:\O_1\to\O_2$ which is a bijection
and the inverse of which $F^{-1}:\O_2\to\O_1$ is also ${\R}cm$
(resp.~${\C}cm$).

A ${\R}cm$ (resp.~${\C}cm$) $n$-manifold over $M_h$
(resp.~$M$) is then a manifold modeled on $M^n_h$ (resp. $M^n$)
with an equivalence class of atlases $(U_1,\var_1)$ such that the
maps $\var_j\var_i^{-1}:\var_i(U_i\cap U_j)\to\var_j(U_i\cap U_j)$
are ${\R}cm$-iso\mor s (resp.~${\C}cm$ iso\mor s)
(See for instance [7] for the definitions in the differentiable
case.)

\bs\noindent{\bf 5.33 \ Remarks.}
a) It is easy to produce examples of ${\R}cm$ or ${\C}cm$
manifolds by gluing open sets in $M^n_h$, or respectively
$M^n$ by ${\R}cm$, or respectively
${\C}cm$ isomorphisms.

b) It seems that some natural candidates for ${\R}cm$ manifolds
don't satisfy the definitions. For instance the unitary group
$U(M)=\{u\in M\mid u^*u=uu^*=I\}$ is such a candidate, but there is a
problem with the charts which one would expect to be connected with the
Cayley transform $u\to i(u-e^{i\th}I)(u-e^{-i\th}I)^{-1}$
or to logarithms. Unfortunately the natural domain for the Cayley
transform are the unitaries such that Ker$(u-e^{-i\th}I)=\{0\}$,
which is not an open set while the range would consist of unbounded
selfadjoint operators affiliated with $M$. Since these charts are the
natural ones, the solution could be to find relaxations of the
definitions so that also these charts be admissible.

\bs{\bf 5.34 \ Cyclicgrad submanifolds.}
To conclude the discussion around cyclomorphy here is another interesting class
of geometric objects that arise in this context: the real (respectively,
complex) cyclicgrad submanifolds in $M^n_h$ (respectively, in $M^n$).
A differentiable (respectively, complex-analytic) submanifold
$V\subset M^n_h$ (respectively, $W\stt M^n$) is real-cyclicgrad
(respectively, complex-cyclicgrad) if for every $Y\in V$
(respectively, $Z\in W$), $T_YV\stt L^2(\d{\Cn}[Y])$
(respectively, for the complex tangent space
$T_ZW\stt L^2(\d{\Cn}[Z])$. Note that these manifolds are perpendicular to
automorphic orbits since the latter have tangent spaces contained in the
orthogonal to cyclic \gra s.

\section{The Lie algebras}

\noindent{\bf 6.1} \ Vect ${\Cn}$ is a Lie \alg \ under the bracket
\[
[P,Q] = (D_PQ_j-D_QP_j)_{1\leq j\leq n}
\]
where $P=(P_j)_{1\leq j\leq n}$, \ $Q=(Q_j)_{1\leq j\leq n}$.
Possible confusion with the commutator on ${\Cn}$ which is denoted the
same way can be avoided by checking the context.

Since the map
\[
{\mb{Vect}}\,{\Cn} \ni P\to D_P\in {\mb{Der}}\,{\Cn}
\]
into the derivations of ${\Cn}$ is injective, the equality
$D_{[P,Q]}=[D_P,D_Q]$ implies that $[P,Q]$ is indeed a Lie \alg \
bracket. The last equality in turn is a straightforward computation
\begin{eqnarray*}
D_PD_Q X_{i_1}\dots X_{i_p} &=& D_P
\sum_{1\leq j\leq p}X_{i_1}\dots X_{i_{j-1}} Q_{i_j}
X_{i_{j+1}}\dots X_{i_p} \\
&=&\sum_{1\leq j,k\leq p\atop j\neq k}
X_{i_1}\dots X_{i_j}P_{ij}
X_{i_{j+1}}\dots X_{i_{k-1}} Q_{i_k}X_{i_k}\dots X_{i_p} \\
&&\qquad + \sum_{1\leq j\leq p}X_{i_1}\dots X_{i_{j-1}}
(D_PQ_{i_j})X_{i_{j+1}}\dots X_{i_p} \ .
\end{eqnarray*}
Combining with the analogues formula for $D_QD_P$ we get
\begin{eqnarray*}
[D_PD_Q] X_{i_1}\dots X_{i_p} &=&
\sum_{1\leq j\leq p}X_{i_1}\dots X_{i_{j-1}}
(D_P Q_{i_j}-D_QP_{i_j})
X_{i_{j+1}}\dots X_{i_p} \\
&=& D_{[P,Q]}
X_{i_1}\dots X_{i_p} \ .
\end{eqnarray*}

\bs\noindent{\bf 6.2} \  We have
$|[P,Q]|_{R,0} \leq |P|_{R,0}|Q|_{R,1}+|P|_{R,1}|Q|_{R,0}$.
In particular if $R'>R$ the definition of the  bracket extends
to $P,Q\in{\mb{Vect}}\,{\C}_{\<n\>,R',0}$ with
$[P,Q]\in {\mb{Vect}}\,{\C}_{\<n\>,R,0}$. This turns
\[
{\mb{Vect}}\,{\C}_{\<n\>,>R}=\bigcup_{R'>R}
{\mb{Vect}}\,{\C}_{\<n\>,R,0}
\]
into a Lie \alg.

\bs\noindent{\bf 6.3 \ Lemma.} {\it Let} \ $H,K\in {\rm{Vect}}\,{\Cn}$ \
{\it and assume}
\begin{eqnarray*}
\Psi_n(H) &\leq & C_1(1-\a X_1)^{-e} \\
\Psi_n(K) &\leq & C_2(1-\a X_1)^{-f} \\
\end{eqnarray*}
Then
\[
\Psi_n([H,K]) \leq C_1C_2 (e+f)(1-\a X_1)^{-e-f-1} \ .
\]

\bs{\bf Proof.} This is a straightforward computation based on
\begin{eqnarray*}
\Psi_n([H,K]) &\leq & \Psi_n((D_{|H|}|K_1|),\dots, (D_{|H|}|K_n|))
   + \Psi_n((D_{|K|}|H_1|),\dots, (D_{|K|}|H_n|))\\
&\leq & D_{\Psi_n(H)} \Psi_n(K) +D_{\Psi_n(K)} \Psi_n(H)\\
&\leq & C_1C_2(D_{(1-\a X_1)^{-e}}(1-\a X_1)^{-f} +
D_{(1-\a X_1)^{-f}}(1-\a X_1)^{-e} )
\end{eqnarray*} \qed

\bs\noindent{\bf 6.4 \ Proposition.} {\it Let} \ $K^{(j)}\in
{\rm{Vect}}\,{\Cn}$,
$0\leq j\leq m$,  \
{\it and let} $M\geq |K^{(j)}|_{R',0}$, $0\leq j\leq m$. {\it If
$R'>R>0$, then}
\[
|{\mb{ad}}\,K^{(n)}{\mb{ad}}\,K^{(n-1)}{\mb{ad}}\,K^{(1)}K^{(0)}|_{R,0}
\ \leq \ M^{m+1}2^mm! (1-R/R')^{-2-1}
\]

\bs{\bf Proof.} The assumption $M\geq |K^{(j)}|_{R',0}$ gives
$\Psi_n(K^{(j)})\leq M(1-X_1/R')^{-1}$. We then prove by induction that
\[
\Psi_n({\mb{ad}}\,K^{(n)}{\mb{ad}}\,K^{(n-1)}\dots
{\mb{ad}}\,K^{(1)}{\mb{ad}}\,K^{(0)})\leq
M^{m+1}2^mm!(1-X_1/R')^{-2m-1}
\]
using Lemma 2.9. Setting $X_1=R$ then yields the desired result. \qed

\bs\noindent{\bf 6.5} \ It is easily seen that the Lie bracket we considered
extends to Vect ${\C}\<\<X_1,\dots ,X_n\>\>$ (see for instance
the discussion fo the grading in 6.9). It is immediate from Lemma 6.3
that the following proposition holds.

\bs{\bf  Proposition.}   {\it The set}
$\{H\in {\mb{Vect}}\,{\C}\<\<X_1,\dots ,X_n\>\>\mid \psi_n(H)\leq
C(1-\a X_1)^{-p}$ for some $C>0$ and $p\in{\N}\}$
{\it is a Lie subalgebra.}

\bs\noindent{\bf 6.6} \  Remark that the Lie \alg \ in Proposition 6.5 is
contained in Vect ${\C}_{\<n\>,R',0}$, where $\a R'<1$, and contain
Vect ${\C}_{\<n\>,R'',0}$, when $\a R''>1$. Combining this with
Proposition 6.5 gives the following result.

\bs{\bf  Corollary.}
Vect\,${\C}_{\<n\>,>R}$ {\it is a Lie sub\alg \ of}
Vect\,${C}\<\<X_1,\dots ,X_n\>\>$.

\bs{\bf 6.7 \ Remark.} Unless in some unexpected way the estimates in
Proposition 6.4 can be vastly improved, it seems unlikely that
convergence of the Campbell-Hausdorff series\linebreak (see [2]) can
be obtained in the operator \alg \ context.

\bs{\bf 6.8} \ We set aside the completions now to point out some simple
features of Vect ${\Cn}$ and we leave  the reader to carry over the
obvious considerations for the completions.

\bs{\bf 6.9}  \ Vect ${\Cn}$ has a grading
$\bigoplus_{k\geq -1}V_k$ where $V_k$ consists of $n$-tuples
$P\!=\!(P_1,\dots ,P_n)$ which are homogeneous of degree $k+1$.
It is easily seen that $[V_k,V_\ell]\stt V_{k+\ell}$.

\bs{\bf 6.10} \ It is easily seen that $V_{-1}$ is isomorphic to the
commutative Lie \alg \ ${\C}^n$. Also $V_{\geq 0}=\bigoplus_{k\geq 0}
V_k$ is a Lie sub\alg \ and $V_{\geq p}=\bigoplus_{k\geq  p} V_k$ where
$p\in{\N}$ are ideals in $V_{\geq 0}$. $V_0$ is also a Lie \alg.
Let $E_{ik}=(\d_{ij}X_k)_{1\leq j\leq n}$ (here $\d_{ij}$ is the
Kronecker symbol). Then $[E_{ab},E_{cd}]=\d_{bc}E_{ad}-\d_{da}E_{cb}$,
i.e., $V_0$ is isomorphic to ${\mathfrak gl}(n,{\C})$ with $E_{ab}$
corresponding to the matrix with $(i,j)$-entry equal to $\d_{ai}\d_{bj}$.
In particular the center of $V_0$ is ${\C}(X_j)_{1\leq j\leq n}$.

It is also easily seen that $V_{\geq 0}$ is the semidirect product of
$V_{\geq 1}$ and $V_0$ ($V_0$ acting on $V_{\geq 1}$).

\bs{\bf 6.11} \ The involution on Vect ${\Cn}$ being defined
component-wise, we have $D_{P^*}Q_j^*=(D_PQ_j)^*$ so that
$P\rightsquigarrow P^*$ is a conjugate-linear auto\mor \ of
Vect ${\Cn}$. Then the selfadjoint part Vect ${\C}^{sa}_{\<n\>}$
is a real Lie-\alg. Remark also that $V_0^{sa}$ is isomorphic to
${\mathfrak gl}(n,{\R})$.

\bs{\bf 6.12}  \ If $\tau:{\Cn}\to{\C}$ is a trace, we define
\[
{\mb{Vect}}\,{\C}_{\<n|\tau\>} =
\{P\in {\mb{Vect}}\,{\Cn}\mid\sum_{1\leq j\leq n}
\tau(P_j(\d_jR))=0 \ {\mb{for all}} \ R\in{\Cn}\} \ .
\]
Then Vect ${\C}_{\<n|\tau\>}$ is a Lie sub\alg \ of Vect ${\Cn}$.
It can be viewed as a noncommutative analogue of the Lie \alg s of
volume-preserving vector fields. The defining relations can be also
written $\tau(D_PR)=0$. If $P,Q\in {\C}_{\<n|\tau\>}$ then
$\tau(D_{[P,Q]}R)=0$ since
\[
\tau(D_{[P,Q]}R)=\tau(D_PD_QR-D_QD_PR) =0 \ .
\]

\bs{\bf 6.13} \ Let $P\in{\Cn}$ and consider
$([P,X_j])_{1\leq j\leq n}\in{\mb{Vect}}\,{\Cn}$
(these brackets are commutators). Then $D_{([P,X_j])_{1\leq j\leq n}}
Q=[P,Q]$, i.e. the derivation defined by this element of
Vect ${\Cn}$ is just the commutator with $P$. We have
\begin{eqnarray*}
[([P,X_1],\dots , [P,X_n]),(P_1,\dots ,P_n)] &= &
(D_{([P,X_j])_{1\leq j\leq n}}P_k)_{1\leq k\leq n}\\
&&\ \ - ([D_{(P_k)_{1\leq k\leq n}}P,X_j])_{1\leq j\leq n}
-([P,P_j])_{1\leq j\leq n} \\
&=& -([D_{(P_k)_{1\leq k\leq n}}P,X_j])_{1\leq j \leq n} \ .
\end{eqnarray*}
So $\{([P,X_j])_{1\leq j\leq n}\in {\mb{Vect}}\,{\Cn}\mid P\in{\Cn}\}$
is an ideal in Vect ${\Cn}$. This ideal corresponds to the inner
derivations. In particular since a trace applied to a commutator gives
zero, we infer that
\[
\{([P,X_j])])_{1\leq j\leq n}\in {\mb{Vect}}\,{\Cn}\mid P\in{\Cn}\}
\stt {\mb{Vect}}\,{\C}_{\<n|\tau\>} \ .
\]

\bs{\bf 6.14} \ Clearly
\[
{\mb{Vect}}\ {\C}^{sa}_{\<n|\tau\>} = {\mb{Vect}}\ {\C}^{sa}_{\<n\>}
\cap {\mb{Vect}}\ {\C}_{\<n|\tau\>}
\]
is a real Lie sub\alg \ of Vect ${\Cn}$. Note also that
\[
\{(i[P,X_j])
_{1\leq j\leq n}\in {\mb{Vect}}\ {\Cn}\mid P\in{\C}^{sa}_{\<n\>}\}\stt
{\mb{Vect}}\ {\C}^{sa}_{\<n|\tau\>}  \ .
\]

\section{The Lie \alg Vect ${\C}_{\<n|\tau\>}$ in the
semicircular case}

{\bf 7.1} \  We return to the context of 1.4 and we will compute
Vect ${\C}_{\<n|\tau\>}$ in the
semicircular case. An important consequence will be that the assumptions
for the exponentiation result, Theorem 3.8, are satisfied in the
semicircular case.

${\mathcal T}({\C}^n)$ identifies on one hand with the $L^2$-space
of $W^*(s_1,\dots ,s_n)$, but also with the
$L^2$-space of the non-selfadjoint \alg \ generated by
$\ell_1,\dots ,\ell_n$, w.r.t.~the scalar product which corresponds to
the vacuum on the $C^*$-\alg \ generated by $\ell_1,\dots ,\ell_n$. Thus
we will be able to look at cyclic gradients w.r.t.~$s_1,\dots ,s_n$ and
$\ell_1,\dots ,\ell_n$ within the same Hilbert space. We shall denote by
\begin{eqnarray*}
&&\d^{\ell} =(\d_1^{\ell},\dots ,\d_n^{\ell})
{\C}\<\ell_1,\dots ,\ell_n\>\to {\mb{Vect}}\,{\C}\<\ell_1,\dots
,\ell_n\>\\
&&{\mb{and}}\\
&&
\d^s=(\d^s_1,\dots ,\d^s_n):{\C}\<s_1,\dots ,s_n\>\to
{\mb{Vect}}\,{\C}\<s_1,\dots ,s_n\>
\end{eqnarray*}
the corresponding cyclic \gra s. Similarly $\p^\ell$ and $\p^s$ will
denote the two free difference quotient \gra s.

\bs\noindent{\bf 7.2 \ Proposition.} \ {\it Assume} $i_j\neq i_{j+1} \
(1\leq j<p)$ {\it and assume} $k_1>0,\dots ,k_p>0$. {\it Then:}
\[
(\p^s_j P_{k_1}(s_{i_1})\dots P_{k_p}(s_{i_p}))(1\ot 1) =
(\p^\ell_j \ell_{i_1}^{k_1}\dots \ell^{k_p}_{i_p})(1\ot 1)
\]
{\it If moreover $i_1\neq i_p$ if $p>1$, then}
\[
(\d^s_j P_{k_1}(s_{i_1})\dots P_{k_p}(s_{i_p}))1=
(\d^\ell_s \ell_{i_1}^{k_1}\dots \ell^{k_p}_{i_p})1 \ .
\]

\bs{\bf Proof.} By the recurrence relation for Gegenbauer \po s which we
did recall in 1.4 we have
\[
(1-rt+r^2)^{-1} =\sum_{n\geq 0} P_n(t)r^n \ .
\]
In particular computing the difference quotient for the left- and
right-hand sides, we get
\begin{eqnarray*}
\sum_{n\geq 0} \frac{P_n(t_1)-P_n(t_2)}{t_1-t_2} \ r^n& = &
r(1-rt_1+r^2)^{-1} (1-rt_2+r^2)^{-1} \\
&=& \sum_{n\geq 1}\left(\sum_{0\leq k\leq n-1}
P_k(t_1)P_{n-1-k}(t_2)\right)r^n  \ .
\end{eqnarray*}
It follows that
\[
\frac{P_n(t_1)-P_n(t_2)}{t_1-t_2}  =\sum_{0\leq k\leq n-1}
P_k(t_1)P_{n-1-k}(t_2)
\]
i.e.,
\[
\p P_n = \sum_{0\leq k\leq n-1} P_k\ot P_{n-1-k} \ .
\]
This in turn gives
\begin{eqnarray*}
&&(\p^s_j P_{k_1}(s_{i_1})\dots P_{k_p}(s_{i_p}))(1\ot 1) \\
&& \ \ \sum_{\{m|i_m=j\}} \ \sum_{0\leq h<k_m}
(P_{k_1}(s_{i_1})\dots P_h(s_{i_m})1)\ot
(P_{k_m-1-h}(s_{i_m})P_{k_{m+1}}(s_{i_{m+1}})\dots
P_{k_p}(s_{i_p})1) \\
&& \ \ =   \sum_{\{m|i_m=j\}} \ \sum_{0\leq h<k_m}
\big( e^{\ot k_1}_{i_1}\ot\dots\ot e^{\ot h}_{i_m}\big)\ot
\big(e_{i_m}^{\ot k_m-1-h}\ot\dots\ot e_{i_p}^{\ot k_p}\big) \\
&& \ \ =  \sum_{\{m|i_m=j\}} \ \sum_{0\leq h<k_m}
\big(\ell^{k_1}_{i_1}\dots \ell^h_{i_m}1\big)\ot
\big(\ell^{k_m-1-h}_{i_m}\dots \ell^{k_p}_{i_p}1\big)\\
&& \ \ =\big( \p^\ell_j \, \ell^{k_1}_{i_1}\dots \ell^{k_p}_{i_p}\big)1\ot 1
\ .\end{eqnarray*}
Similarly under the additional assumption that
$i_p\neq i_1$ in case $p>1$, we have
\begin{eqnarray*}
&& \d^s_j P_{k_1}(s_{i_1})\dots P_{k_p}(s_{i_p})1 \\
&& \ = \sum_{\{m|i_m=j\}} \ \sum_{0\leq h<k_m}
P_{k_m-1-h}(s_{i_m})P_{k_m+1}(s_{i_{m+1}})\dots
P_{k_p}(s_{i_p})P_{k_1}(s_{i_1})\dots
P_{k_{m-1}}(s_{i_{m-1}})P_h(s_{i_m})1\\
&& \ = \sum_{\{m|i_m=j\}} \ \sum_{0\leq h<k_m}
e_{i_m}^{\ot k_m-1-h}\ot e_{i_{m+1}}^{\ot k_{m+1}}\ot\dots\ot
e_{i_p}^{\ot k_p} \ot e_{i_1}^{\ot k_1}\ot\dots \ot
e_{i_{m-1}}^{\ot k_{m-1}}\ot e_{i_m}^{\ot h}\\
&& \ =  \sum_{\{m|i_m=j\}}
\ell_{i_m}^{k_m-1-h}\ell_{i_{m+1}}^{k_{m+1}}\dots
\ell_{i_p}^{k_p}\ell_{i_1}^{k_1}\dots \ell_{i_{m-1}}^{k_{m-1}}
\ell^h_{i_m}1 \\
&& \ = \big(\d^\ell_j\, \ell_{i_1}^{k_1}\dots
  \ell_{i_p}^{k_p}\big)1 \ .
\end{eqnarray*} \qed

\bs\noindent{\bf 7.3 \ Lemma.} {\it Let}
${\F}_k\stt{\C}\<s_1,\dots ,s_n\>$ {\it be the linear span of
$P_{k_1}(s_{i_1})\dots P_{k_p}(s_{i_p})$ with $k_j>0$
$(1\leq j\leq p)$, $k_1+\dots +k_p=k$, $i_j\neq i_{j+1}$
$(1\leq j< p)$ and $i_1\neq i_p$ \ if \ $p>1$, and let
${\F}=\sum_{k\geq 0}{\F}_k$. Then} \
${\F}+$ Ker\, $\d^s={\C}\<s_1,\dots ,s_n\>$.

{\bf Proof.} Recall from 1.3 that Ker $\d^s={\mb{Ker}}\, C^s$ where
$C^s$ is the cyclic symmetrization. Moreover Ker $\d^s=\sum_{k\geq 0}
{\mb{Ker}}\,\d^{s,k}$ is a direct sum decomposition where
$\d^{s,k}$ is the restriction of $\d^s$ to the homogeneous \ncc \po s of
degree $k$ and similarly
\[
{\mb{Ker}}\,C^s=\sum_{k\geq 0} {\mb{Ker}}\,C^{s,k}
\]
so that \ Ker $\d^{s,k}={\mb{Ker}}\,C^{s,k}$.

Let ${\G}_k$ be the linear span of
$s^{k_1}_{i_1}\dots s^{k_p}_{i_p}$ with $k_j>0$
$(1\leq j\leq p)$, $k_1+\dots +k_p=k$, $i_j\neq i_{j+1}$
$(1\leq j< p)$ and $i_1\neq i_p$ \ if \ $p\geq 2$, and put
${\G}=\sum_{k\geq 0}{\G}_k$. It is easily seen that
\[
{\G}_k + {\mb{Ker}}\,C^{s,k}=({\C}\<s_1,\dots ,s_n\>)_k \ .
\]
It will suffice to prove that
\[
\sum_{0\leq k\leq m} ({\F}_k+{\mb{Ker}}\,C^{s,k})  \ = \
\sum_{0\leq k\leq m} ({\G}_k+{\mb{Ker}}\,C^{s,k}) \ .
\]
This is shown by induction on $m$, using
\[
{\F}_m+\sum_{0\leq k< m} ({\C}\<s_1,\dots ,s_n\>)_k \ = \
{\G}_m+\sum_{0\leq k< m} ({\C}\<s_1,\dots ,s_n\>)_k
\]
which follows from $P_k(s)=cs^k$ (mod $\sum_{0\leq h < k}({\C}\<s\>)_h$) for
some
$c\neq 0$.\qed

\bs\noindent{\bf 7.4 Theorem.}
$(\d^s{\C}\<s_1,\dots ,s_n\>)(1\oplus\dots\oplus 1) =
(\d^\ell {\C}\<\ell_1,\dots ,\ell_m\>)(1\oplus\dots\oplus 1)$.

\bs{\bf Proof.} Let ${\mathcal L}_k$ be the linear span of
$\ell^{k_1}_{i_1}\dots \ell^{k_p}_{i_p}$, $k=k_1+\dots +k_p$,
$i_j\neq i_{j+1}$, $(1\leq j\leq p)$, $i_1\neq i_p$ if $p\geq 2$,
and ${\LL}=\sum_{k\geq 0}{\LL}_k$. Then
${\LL}+{\mb{Ker}}\,C^\ell={\C}\<\ell_1,\dots ,\ell_n\>$ and hence
${\LL}+{\mb{Ker}}\,\d^\ell={\C}\<\ell_1,\dots ,\ell_n\>$.

It follows that
\begin{eqnarray*}
&&(\d^\ell{\C}\<\ell_1,\dots ,\ell_n\>)(1\oplus\dots\oplus 1)\\
&& \ = (\d^\ell{\LL})(1\oplus\dots\oplus 1)\\
&& \ = (\d^s{\F})(1\oplus\dots\oplus 1) \\
&& \ = (\d^s {\C}\<s_1,\dots ,s_n\>)(1\oplus\dots\oplus 1)
\end{eqnarray*}
where we used 7.2 and 7.3 in the last two steps. \qed

\bs\noindent{\bf 7.5 \ Theorem.} {\it Identifying
${\T}({\C}^n)$ and $L^2(W(s_1,\dots ,s_n),\tau)$ we have:}
\[
({\T}({\C}^n))^n \ominus
\overline{\d^s{\C}\<s_1,\dots ,s_n\>} \ = \
\{((\ell^*_j-r^*_j)\xi)_{1\leq j\leq n}\mid
\xi\in {\T}({\C}^n)\} \ .
\]
{\it In particular}
\[
{\mb{Vect}}\,{\C}\<s_1,\dots ,s_n\mid\tau\> =\sum_{k\geq 0}
{\mathcal X}_k
\]
{\it where ${\mathcal X}_k\stt ({\T}_k({\C}^n))^n$, where
${\mathcal X}_0=0$ and ${\mathcal X}_k$ for $k\geq 1$ is spanned by}
\[
\big(\d_{j,i_0}e_{i_1}\ot\dots\ot e_{i_k}-\d_{j,i_k}
e_{i_0}\ot\dots\ot e_{i_{k-1}}\big)_{1\leq j\leq n}
\]
{\it where} \ $(i_0,\dots ,i_k)\in \{1,\dots ,n\}^{k+1}$.

\bs{\bf Proof.} $\d^s{\C}\<s_1,\dots ,s_n\>$ identifies with
\[
\d^s{\C}\<s_1,\dots ,s_n\>)(1\oplus\dots\oplus 1) =
(\d^\ell{\C}\<\ell_1,\dots ,\ell_n\>)(1\oplus\dots\oplus 1) \ .
\]
To compute \ $({\T}({\C}^n))^n \ominus
\overline{\d^\ell{\C}\<\ell_1,\dots ,\ell_n\>(1\oplus\dots\oplus 1)}$
we use the exact sequence 1.3. The maps in the exact sequence are
homogeneous, so we get exact sequences
\[
({\Cn})_k  \ \stk{\d}{\longrightarrow} \ V_{k-2}
\ \stk{\th}{\longrightarrow} \ ({\Cn})_k
\]
(recall that $V_k$ are the \ncc vector fields with components
homogeneous of degree $k+1$). Endowing $({\Cn})_k$ with the scalar
product in which the monomials $X_{i_1}\dots X_{i_k}$ form an orthonormal
basis we have isometries
\begin{eqnarray*}
&& ({\Cn})_k  \ \longrightarrow \ {\T}_k({\C}^n)\\
&& V_{k-2} \ \longrightarrow \ ({\T}_{k-1}({\C}^n))^n
\end{eqnarray*}
which map $X_{i_1}\dots X_{i_k}$ to $e_{i_1}\ot\dots\ot e_{i_k}$.
In this correspondence $\d$ identifies with $\d^\ell$ acting in
the $L^2$-space and $\th$ with the map
\[
\th^{\ell}((\xi_j)_{1\leq j\leq n})=\sum_{1\leq j\leq n}
(\ell_j-r_j)\xi_j \ .
\]
Clearly, then
\[
({\T}_k({\C}^n))^n\ominus \d^\ell
({\C}\<\ell_1,\dots ,\ell_n\>)_{k+1} =
(\th^\ell)^*\xi = ((\ell^*_j-r^*_j)\xi)_{1\leq j\leq n} \ .
\]
Taking $\xi=e_{i_0}\ot e_{i_1}\ot\dots\ot e_{i_k}$ we get the last part
of the statement. \qed

\bs{\bf 7.6.} Transcribing the last part
of the preceding theorem in terms of the $s_j$'s instead of the
$e_j$'s, we get the following result.

\bs{\bf  Corollary} \ {\it The elements}
\[
F_I=\big(\d_{i_0,j}P_{k_0-1}(s_{i_0})P_{k_1}(s_{i_1})\dots
P_{k_p}(s_{i_p})-
\d_{i_p,j}P_{k_0}(s_{i_0})P_{k_1}(s_{i_1})\dots
P_{k_p-1}(s_{i_p})\big)_{1\leq j\leq n}
\]
{\it where} $I=(\underbrace{i_0,\dots ,i_0}_{k_0{\mb{\footnotesize
-times}}} ,\dots ,\underbrace{i_p,\dots ,i_p}_{k_p{\mb{\footnotesize
-times}}})$,
$k_r>0$ $(0\leq r\leq p)$, {\it span} Vect\,${\C}\<s_1,\dots
,s_n\mid\tau\>$.

\bs{\bf 7.7.} We will now clarify the relations among the elements
$F_I$. On ${\T}_k({\C}^n)$ let $R_k$ be the cyclic permutation
$R_k e_{i_1}\ot\dots \ot e_{i_k}= e_{i_k}\ot e_{i_1}\ot\dots \ot
e_{i_{k-1}}$ (in particular $R_01=0$, $R_1e_j=e_j$ and
$R_2e_i\ot e_j=e_j\ot e_i$). Let
$R=\bigoplus_{k\geq 1} R_k$ be the  operator on
${\T}({\C}^n)$. Let also $P$ be the projection onto ${\C}1$.

\bs{\bf  Lemma} \ {\it Let \ $T_j=\ell_j-r_j$. Then}
$\ds{\sum_{1\leq j\leq n}} T_jT^*_j = 2I-2P-(R+R^*)$.
{\it In particular}
\[
{\mb{Ker}}\,(\th^\ell)^* = {\mb{Ker}}\sum_{1\leq j\leq n} T_jT^*_j =
{\mb{Ker}}(I-P-R) \ .
\]

\medskip{\bf Proof.}
\begin{eqnarray*}
\sum_{1\leq j\leq n} T_jT^*_j &=&
\sum_{1\leq j\leq n} \big(\ell_j\ell^*_j +r_jr^*_j -
\ell_jr_j^* -r_j\ell^*_j\big) \\
&=& (I-P)+(I-P)-R-R^* \ .
\end{eqnarray*}
Since $R\left.\right|_{{\T}_+}$ is a unitary operator, we have\\
Ker$((2I-R-R^*)|_{{\T}_+} ={\mb{Ker}}(I-R)|_{{\T}_+}$. \qed

\bs{\bf 7.8. Remark.} Let $N$ be the number operator, i.e.
$N1=0$, $Ne_{i_1}\ot\dots\ot e_{i_k}=ke_{i_1}\ot\dots\ot e_{i_k}$. Then
\[
{\mb{Ker}}((I-R)|_{{\T}_+} ={\mb{Ker}}(N-C)|_{{\T}_+}
\]

\bs{\bf 7.9. Remark.} To get a basis of Vect\,${\C}\<s_1,\dots ,
s_n|\tau\>$ we shall use a basis of \linebreak $(I-R)({\T}_+({\C}^n))$.
Consider on $\{1,\dots ,n\}^{k+1}$ the lexicographic order
$\prec$ and let
\[
\O(k+1)\!=\!\{(i_0,\dots ,i_k)\!\in\! \{i,\dots ,n\}^{k+1}\!\mid\!
(i_0,\dots ,i_k)\!\prec \!(i_k,i_0,\dots ,i_{k-1}) \ {\mb{and}} \
(i_0,\dots ,i_k)\neq\! (i_k,i_0,\dots ,i_{k-1})\}
\]
Then
\[
\Psi=\bigcup_{k\geq 1}
\{F_{i_0\dots i_k}-F_{i_ki_0\dots i_{k-1}}\mid
(i_0\dots i_k)\in\O(k+1)\}
\]
is a basis of Vect\,${\C}\<s_1,\dots ,s_n|\tau\>$ over ${\C}$.
Indeed the set
\[
\Phi=\bigcup_{k\geq 1}
\{e_{i_0}\ot\dots\ot e_{i_k}-e_{i_k}\ot e_{i_0}\dots e_{i_{k-1}}\mid
(i_0\dots i_k)\in\O(k+1)\}
\]
is a homogeneous basis of $(\bigoplus_{k\geq 2} {\T}_k({\C}^n)\ominus
{\mb{Ker}}(I-P-R)$ which in view of Lemma 7.7 is mapped by
$(\th^\ell)^*$ to Vect\,${\C}\<s_1,\dots ,s_n|\tau\>$.
We have used here that $(\th^\ell)^*\xi=(T_j^*\xi)_{1\leq j\leq n}$,
so that Ker$(\th^\ell)^* ={\mb{Ker}}\,\sum_{1\leq j\leq n} T_jT^*_j$.
Also $F_{i_0\dots i_k}=(\th^\ell)^*e_{i_0}\ot\dots\ot e_{i_k}$.

\bs{\bf 7.10. Remark.} The same kind of considerations as in 7.9
yield also another natural basis. If $I\!=\!(i_0,\dots ,i_k)\in\{1,\dots
,n\}^{k+1}$, let per$(I)$ be the least non-zero cyclic period of~$I$,
i.e., the least $m\!\in\!\{1,\dots ,k\}$ such that $i_s=i_t$ whenever
$s\!\equiv\! t$ (mod $m$). If $\zeta^m=1$,\linebreak $\zeta\neq 1$ and
$m={\mb{per}}(I)$ let $F(I,\zeta)=\sum_{0\leq j<m} \zeta^j
F_{(i_j,i_{j+1},\dots ,i_k,i_0,\dots ,i_{j-1})}$. Let also
$\o(k+1)=\{(i_0\dots i_k)\in \{1,\dots ,h\}^{k+1}\mid
(i_0,\dots i_k)\prec (i_j,i_{j+1},\dots ,i_k,i_0,\dots ,i_{j-1}),
j=1,\dots ,k\}$ and $\rho(m)=\{\zeta\in{\C}\mid\zeta\neq 1, \
\zeta^m=1\}$. Then
\[
\bigcup_{k\geq 1}\{F(I,\zeta)\mid I\in \o(k+1), \
\zeta\in\rho ({\mb{per}}(I))\}
\]
is a basis of Vect\,${\C}\<s_1,\dots ,s_n|\tau\>$.

\bs{\bf 7.11.} To construct a basis over ${\R}$ of the
real Lie \alg \ Vect\,${\C}^{sa}\<s_1,\dots ,s_n|\tau\>$
we can take the hermitian and antihermitian parts of the elements of the
previous bases. In doing so one must pay attention to the fact that
if per$(I)=2$ then $F_{i_0\dots i_k}-F_{i_ki_0\dots i_{k-1}}$
is antihermitian. For instance the basis in 7.9 then gives the basis:
\begin{eqnarray*}
&&\bigcup_{k\geq 1}\left\{  F_{i_0\dots i_k}+F_{i_k\dots i_0}-
F_{i_ki_0\dots i_{k-1}}-F_{{i_{k-1}}\dots i_0i_{k-1}}\mid
  (i_0,\dots ,i_k)\in \O(k+1),
{\mb{per}}(i_0,\dots ,i_k) >2  \right\}\\
&&\cup
\bigcup_{k\geq 1}\left\{ \sqrt{-1} (F_{i_0\dots i_k}-F_{i_k\dots i_0}
-F_{i_ki_0\dots i_{k-1}}-F_{{i_{k-1}}\dots i_0i_{k-1}}\mid
  (i_0,\dots ,i_k)\in \O(k+1)\right\}
\end{eqnarray*}

\bs{\bf 7.12.} Putting together the results of this section and the
results about exponentiation and automorphic orbits we see that we
have proved the following result.

\bs{\bf Theorem} \ {\it If $s_1,\dots ,s_n$ is a semicircular system
and $M=W^*(s_1,\dots ,s_n)$ then}\\
$T\!AO_\i(s_1,\dots ,s_n)\cap{\mb{Vect}}\,{\C}^{sa}\<s_1,\dots ,s_n\> =
{\mb{Vect}}\,{\C}^{sa}\<s_1,\dots ,s_n\mid\tau\>$
{\it is dense in $(\d {\C}^\perp_{\<n\>})(s_1,\dots ,s_n)$.
In particular for $K$ in this dense subset of
$(\d {\C}^\perp_{\<n\>})(s_1,\dots ,s_n)$, $D^\e_K$ exponentiates to
a one-parameter auto\mor \ group of $W^*(s_1,\dots ,s_n)$.
Moreover $C^*(s_1,\dots ,s_n)$ is invariant under these}
$\exp(tD^\e_K)$.

\bs{\bf 7.13. Remark} \ It is easily seen that
$V_0\cap {\mb{Vect}}\ {\C}^{sa}\<s_1,\dots ,s_n\mid\tau\>=0$
and\linebreak $V_1\cap {\mb{Vect}}\ {\C}^{sa}\<s_1,\dots ,s_n\mid\tau\>$
is isomorphic to the Lie algebra of the orthogonal group
${\mathfrak o}(n,{\mathbb R})$.

\section{$B$-morphic maps}

{\bf 8.1} \ \ Cyclomorphic maps have a natural generalization to the context of
$B$-von Neumann \alg s, i.e.~the complex field is replaced by $B$.
More precisely we deal with $(M,\tau)$ von Neumann \alg s with specified normal
faithful trace-state and there is a specified unital inclusion
$B\hookrightarrow M$ so that $B$ is a von Neumann sub\alg . For lack of better
names the analogues of cyclomorphic maps will be called $B$-morphic.
If $(Y_j)_{1\leq j\leq n}\in M^n_h$ we shall denote by
$\D_h(Y_1,\dots ,Y_m:B_h)$ the norm-closed linear span in
${\mathcal L}(M^n_h,B_h)$ of the maps
\begin{eqnarray*}
M^n_h  \ni m &= &(m_j)_{1\leq j\leq n}\rightsquigarrow
T(m)+(T(m))^* \\
{\mb{ and}} &&\\
M^n_h \ni m& =& (m_j)_{1\leq j\leq n} \rightsquigarrow
\sqrt{-1} \ (T(m)-(T(m))^*
\end{eqnarray*}
where $T(m)$ is of the form
\[
  T(m)  = \sum_{1\leq j\leq p} E_B(b_0
Y_{i_1}b_1\dots Y_{i_{j-1}}b_{j-1} m_{i_j} b_j Y_{i_{j+1}}b_{j+1}\dots
Y_{i_p}b_p) \ .
\]

Similarly when $(Y_j)_{1\leq j\leq n}\in M^n$ (no selfadjointness requirement)
then $\Delta(Y_1,\dots ,Y_n:B)$ denotes the norm-closed linear span
in ${\mathcal L}(M^n,B)$ of the maps $m\rightsquigarrow T(m)$. Note that
$T(\dt)$ is the differential at $(Y_1,\dots ,Y_n)$ of the map
$M^n\ni (X_j)_{1\leq j\leq n}\rightsquigarrow E_B(b_0X_{i_1}b_1\dots
X_{i_n}b_n)$. Remark that $\Delta(Y_1,\dots ,Y_n:B)$ is a $B\!-\!B$-bimodule,
while $\Delta_h(Y_1,\dots ,Y_n:B_h)$ is stable under $L\rightsquigarrow
bL+Lb^*$, \ $L\rightsquigarrow (bL-Lb^*)\sqrt{-1}$. This last feature is why we
will prefer to work with $\Delta_h(Y_1,\dots ,Y_n:B)$ in the selfadjoint case,
which is the norm closure of the linear span of the $T(\dt)$ in
${\mathcal L}(M^n,B)$.

By $\Delta^\perp_h(Y_1,\dots ,Y_n:B)$ we denote the set of
$K=(K_1,\dots ,K_n)\in M^n_h$ so that $T(K)=0$ for all $T\in \Delta_h(Y_1,\dots
,Y_n:B)$. Similarly, $\Delta^\perp(Y_1,\dots ,Y_n:B)$ is the set of $K\in
M^n$ so that $T(K)=0$ for all $T\in \Delta(Y_1,\dots ,Y_n:B)$.

Our discussion of $B$-morphic maps will be rather sketchy in this paper, we
just want to explain the main idea of the generalization.

\bs{\bf 8.2 \ Definition.} Let $\O\stt M^n_h$ be an open set.
A differentiable map $f:\O\to B$ is hermitian $B$-morphic (abbreviated
$hBm$) if for every $Y\in\O$, $df[Y]\in \Delta_h(Y_1,\dots,Y_n:B)$.
A differentiable map $f:\O\to B$ is weakly hermitian $B$-morphic (abbreviated
$whBm$) if $(df[Y])(\D^\perp_h(Y_1,\dots,Y_n:B))=0$.

\bs{\bf 8.3 \ Definition.} Let $\O\stt M^n_h$ be an open set.
A differentiable map $F:\O\to N^p_h$ where $N$ is also a
$B$-von Neumann \alg \ is $hBm$ if for every $hBm$ map
$f:\o\to B$, $\o\stt N^p_h$ open, the map $f\circ (F\mid_{F^{-1}(\o)\cap \O})$
is $hBm$. The differentiable map $F:\O\to N^p_h$ is $whBm$ if
\[
(dF[Y])(\D^\perp_h(Y:B))\stt \D^\perp_h(F(Y):B) \ .
\]

\bs{\bf 8.4 \ Definition.} Let $\O\stt M^n$ be an open set.
A complex analytic map $f:\O\to B$ is $B$-morphic (abbreviated $Bm$)
if for each $Y\in\O$ the complex differential is such that
$dF[Y]\in \D(Y_1,\dots,Y_n:B)$. A complex analytic map $f:\O\to B$ is weakly
$B$-morphic (abbreviated $wBm$) if for each $Y\!\in\!\O$ the
  complex differential is such that $dF[Y])(\D^\perp_h(Y_1,\dots ,Y_n:B))=~0$.

\bs{\bf 8.5 \ Definition.} Let $\O\stt M^n$ be open and let
$F:\O\to N^p$ be a  complex analytic map. $F$ is $Bm$ if for every $Bm$ map
$f:\o\to B$, $\o\stt
N^p$ open, the map $f\circ F\mid_{F^{-1}(\o)\cap \O}$
is $Bm$.  The map $F$ is $wBm$ if
\[
dF[Y](\D^\perp (Y:B))\stt \D^\perp(F(Y):B) \ .
\]

\bs{\bf 8.6} \ \ If $X_1,\dots ,X_n$ are indeterminates the
\ncc \po s with coefficients in $B$, i.e. $B*_{\C}{\C}\<X_1,\dots ,X_n\>$ will
be denoted by $B\<X_1,\dots ,X_n\>$ or $B_{\<n\>}$.
$B_{\<n\>}$ is spanned by \ncc monomials $b_0X_{i_1}b_1\dots X_{i_p}b_p$.
The involution is
defined so that $(b_0X_{i_1}b_1\dots X_{i_p}b_p)^*=
b^*_pX_{i_p}\dots b_1^* X_{i_1}b^*_0$.

\bs{\bf 8.7 \ Examples.} a) If \ $P\in B_{\<n\>}$ then
$f_P:M^n_h\to B$ defined by
$f_P(Y) = E_B(P(Y))$ is $hBm$.

b) The $hBm$ maps of $\O$ to $B$ form an  algebra.

c) If \ $P\in B_{\<n\>}$ then $g_P:M^n\to B$ defined by
$g_P(Z)=E_B(P(Z))$ is a $Bm$ map.

d) The $Bm$ maps $\O\to B$ form an algebra.

e) If $P\!=\!(P_j)_{1\leq j\leq m}\in (B_{\<n\>})^m$, then
$G_P:M^n\to M^m$ defined by $G_P(Z)=(P_j(Z))_{1\leq j\leq m}$ is a
$Bm$ map.

f) If $P=(P_j)_{1\leq j\leq m}\in (B_{\<n\>})^n$ and
$P_j=P^*_j$, $1\leq j\leq m$, then
$F_P:M^n_h\to M^m_h$ defined by $F_P(Y)=(P_j(Y))_{1\leq j\leq m}$
is a $hBm$ map.

\bs{\bf 8.8 \ Theorem.} {\it Assume we have
$I\in B\stt N\stt M$ and let $Y_1,\dots ,Y_m\in M_h$ be such that
$\{Y_1,\dots ,Y_m\}$  and $N$ are freely independent over $B$
in $(M,E_B)$. Then the map $\Phi:N^n_h\to M_h^{n+m}$ defined by
$\Phi(T_1,\dots ,T_n)=(T_1,\dots ,T_n,Y_1,\dots ,Y_m)$ is a
$hBm$ map.}

\bs{\bf 8.9 \ Theorem.} {\it Assume we have
$I\in B\stt N\stt M$ and let $Z_1,\dots ,Z_m\in M$ be such that
$\{Z_1,\dots ,Z_m\}$  and $N$ are freely independent over $B$
in $(M,E_B)$. Then the map $\Psi:N^n\to M^{n+m}$ defined by
$\Psi(T_1,\dots ,T_n)=(T_1,\dots ,T_n,Z_1,\dots ,Z_m)$ is a
$Bm$ map.}

\newpage
\begin{center}{\bf References}\end{center}
\begin{description}
\item{[1]}
P.Biane and D.Voiculescu, A free \pr ty analogue of the Wasserstein
metric on trace-state space. Preprint, Berkeley (2000).

\item{[2]} N.Bourbaki, Groupes et alg\`ebres de Lie, chapitres II et III.
(Hermann, Paris, 1972).

\item{[3]}
O.Bratteli and D.W.Robinson, {\it Operator Algebras and Quantum
Statistical Mechanics, I}. (Springer-Verlag, 1979).

\item{[4]}
T.Cabanal-Duvillard and A.Guionnet, Large deviations, uper bounds and \ncc
entropies for some matrices ensembles. Preprint, Paris (1999).

\item{[5]} T.Cabanal-Duvillard and A.Guionnet,
Around Voiculescu's free entropies. Preprint, ENS Lyon (2001).

\item{[6]} E.Hille and R.Phillips, {\it Functional Analysis and Semi-Groups},
revised edition (AMS Collqu. Publ., vol.~{\bf 31}, Providence, RI, 1957).

\item{[7]} S.Lang, {\it Introduction to Differentiable Manifolds},
(1962).

\item{[8]} G.-C.Rota, B.Sagan, and P.R.Stein, A cyclic derivative in
\ncc \alg . {\it J. of Algebra} {\bf 64} (1980), 54--75.

\item{[9]}
N.Ya.Vilenkin, {\it Special Functions and the Theory of Group
Representations}. (Moscow, 1965).

\item{[10]}
D.Voiculescu, Symmetries of some reduced free product $C^*$-algebras, in
{\it Operator Algebras and Their Connections with Topology and
Ergodic Theory}, Lecture Notes in Math., vol.~{\bf 1132},
(Springer, 1985), 556--588.

\item{[11]} D.Voiculescu, Circular and semicircular systems and
free product factors, in {\it Operator Algebras, Unitary Representations
and Invariant Theory}, Progress in Mathematics {\bf 92}, (Birkhauser, 1990),
45--60.

\item{[12]} D.Voiculescu, The analogues of entropy and of Fisher's information
measure in free \pr ty theory, V: Noncommutative Hilbert transforms.
{\it Invent. Math.} {\bf 132} (1998), 182--227.

\item{[13]} D.Voiculescu, A note on cyclic \gra s.
{\it Indiana Univ. Math. J.} {\bf 49}, no.~3 (2000), 837--841.

\item{[14]} D.Voiculescu, Free entropy. Preprint, Berkeley (2001).

\item{[15]} D.Voiculescu, K.J.Dykema and A.Nica, {\it Free Random
Variables}, CRM Monograph Series vol.~1. (AMS, Providence, RI, 1992).
\end{description}
\end{document}